\documentstyle[amsthm]{article}
\newtheorem{th}{Theorem}
\newtheorem{lm}{Lemma}

\newtheorem{exs}{Examples}
\newtheorem{dfn}{Definition}
\newtheorem{re}{Remark}

\begin{document}
\author{R.Campoamor-Stursberg \and V.O.Manturov}
\title{Invariant Tensors Formulae via Chord Diagrams}
\date{}

\maketitle

\begin{abstract}
We provide an explicit algorithm  to calculate invariant tensors for the adjoint representation of the simple Lie algebra $sl(n)$, as well as arbitrary representation in terms of roots. We also obtain explicit formulae for the adjoint representations of the orthogonal and symplectic Lie algebras $so(n)$ and $sp(n)$.
\end{abstract}

\section{Introduction}

In the last few years, various knot invariants have been
discovered. Nowadays, Vassiliev knot invariants are the strongest
among all known ones, see e.g. \cite{Ma2}. The essential
notion in Vassiliev's theory is the notion of chord diagram.

Chord diagrams have interesting algebraic structure.
It turns out \cite{BN} that this structure is deeply connected
with Lie algebra representations.

In his article \cite{BN}, D.Bar--Natan shows the way for constructing
invariant tensors for the coadjoint action on semisimple Lie algebra
tensors by using chord diagrams.

The article \cite{CV} describes the case of an arbitrary representation of the algebra
$sl(2)$, and in \cite{BN2} the connection between this construction and the
adjoint representation of $sl(2)$ is shown, in terms of which the Four
Colour problem can be reformulated in tensorial language; the main features for the
case of the adjoint representation for algebras $sl(n)$ are also described.

The aim of this work is to study the connection between the algebras of Chinese and chord diagrams
and the invariant tensors algebra and finding explicit formulae.
Note, that chord diagrams describe the Vassiliev knot invariants. This provides a connection between knot and representation theory.

Moreover Chinese character (and chord) diagrams describe
the Vassiliev finite--type invariants of knots.
Thus, this theory connects knot theory with representation theory.

The method, proposed by Bar--Natan, requires step--by--step
contraction of concrete tensors. Each step is a contraction
and can be described by a concrete formula.

The main results of the present paper are:

\begin{enumerate}

\item To provide an explicit algorithm for calculating invariant
tensors for the adjoint representation of $sl(n)$ (see also \cite{BN2}).
The result of the algorithm action is a scalar function
of the variable $n$. It turns out that this function is a polynomial of $n$ whose
coefficients are functions on chord diagrams, invariant under the so--called 4T--relation.

One gets some properties of the initial chord diagram from its polynomial.
So, the highest possible power
of polynomials is obtained only for so--called $d$--diagrams,
see \cite{Ma}.

\item Formulae for all representations of $sl(n)$ in the
terms of roots.

\item Explicit formulae for the adjoint representation of
$so(n)$ and for the $k+1$--dimensional representation of $A_1$.

\item Explicit formulae for the adjoint representation of
Lie algebras $sp(n)$.
\end{enumerate}

A shortened version of this paper not including results
on $sp(n)$ is published in \cite{MaVi}.

\begin{dfn}
A {\em chord diagram} (CD) \cite{BN} is a graph that
consists of an oriented cycle (also called {\em circle})
and nonoriented edges, connecting points belonging
to circle (also called {\em chords}). Each graph
vertex is incident to just one circle;
chord diagrams are considered as combinatorial
objects (i.e. up to graph isotopy preserving the
circle orientation).

\end{dfn}

The {\em degree} of a chord diagram is the number of
its chords.

\begin{figure}
\unitlength 0.7mm
\begin{center}
\begin{picture}(155,35)
\multiput(20,35)(0.99,-0.10){3}{\line(1,0){0.99}}
\multiput(22.97,34.70)(0.36,-0.11){8}{\line(1,0){0.36}}
\multiput(25.83,33.82)(0.22,-0.12){12}{\line(1,0){0.22}}
\multiput(28.45,32.39)(0.13,-0.11){17}{\line(1,0){0.13}}
\multiput(30.74,30.47)(0.12,-0.15){16}{\line(0,-1){0.15}}
\multiput(32.60,28.14)(0.11,-0.22){12}{\line(0,-1){0.22}}
\multiput(33.96,25.48)(0.12,-0.41){7}{\line(0,-1){0.41}}
\multiput(34.77,22.60)(0.11,-1.49){2}{\line(0,-1){1.49}}
\multiput(35,19.63)(-0.09,-0.74){4}{\line(0,-1){0.74}}
\multiput(34.62,16.66)(-0.12,-0.35){8}{\line(0,-1){0.35}}
\multiput(33.67,13.83)(-0.11,-0.20){13}{\line(0,-1){0.20}}
\multiput(32.18,11.24)(-0.12,-0.13){17}{\line(0,-1){0.13}}
\multiput(30.20,9)(-0.15,-0.11){16}{\line(-1,0){0.15}}
\multiput(27.82,7.20)(-0.24,-0.12){11}{\line(-1,0){0.24}}
\multiput(25.13,5.90)(-0.41,-0.11){7}{\line(-1,0){0.41}}
\multiput(22.24,5.17)(-1.49,-0.07){2}{\line(-1,0){1.49}}
\multiput(19.25,5.02)(-0.74,0.11){4}{\line(-1,0){0.74}}
\multiput(16.30,5.46)(-0.31,0.11){9}{\line(-1,0){0.31}}
\multiput(13.49,6.49)(-0.20,0.12){13}{\line(-1,0){0.20}}
\multiput(10.94,8.04)(-0.13,0.12){17}{\line(-1,0){0.13}}
\multiput(8.75,10.07)(-0.12,0.16){15}{\line(0,1){0.16}}
\multiput(7.01,12.50)(-0.11,0.25){11}{\line(0,1){0.25}}
\multiput(5.78,15.22)(-0.11,0.49){6}{\line(0,1){0.49}}
\put(5.12,18.13){\line(0,1){2.99}}
\multiput(5.04,21.12)(0.10,0.59){5}{\line(0,1){0.59}}
\multiput(5.56,24.06)(0.11,0.28){10}{\line(0,1){0.28}}
\multiput(6.65,26.84)(0.12,0.18){14}{\line(0,1){0.18}}
\multiput(8.27,29.35)(0.12,0.12){18}{\line(0,1){0.12}}
\multiput(10.36,31.49)(0.16,0.11){15}{\line(1,0){0.16}}
\multiput(12.83,33.17)(0.28,0.12){10}{\line(1,0){0.28}}
\multiput(15.58,34.33)(0.74,0.11){6}{\line(1,0){0.74}}
\multiput(60,35)(0.99,-0.10){3}{\line(1,0){0.99}}
\multiput(62.97,34.70)(0.36,-0.11){8}{\line(1,0){0.36}}
\multiput(65.83,33.82)(0.22,-0.12){12}{\line(1,0){0.22}}
\multiput(68.45,32.39)(0.13,-0.11){17}{\line(1,0){0.13}}
\multiput(70.74,30.47)(0.12,-0.15){16}{\line(0,-1){0.15}}
\multiput(72.60,28.14)(0.11,-0.22){12}{\line(0,-1){0.22}}
\multiput(73.96,25.48)(0.12,-0.41){7}{\line(0,-1){0.41}}
\multiput(74.77,22.60)(0.11,-1.49){2}{\line(0,-1){1.49}}
\multiput(75,19.63)(-0.09,-0.74){4}{\line(0,-1){0.74}}
\multiput(74.62,16.66)(-0.12,-0.35){8}{\line(0,-1){0.35}}
\multiput(73.67,13.83)(-0.11,-0.20){13}{\line(0,-1){0.20}}
\multiput(72.18,11.24)(-0.12,-0.13){17}{\line(0,-1){0.13}}
\multiput(70.20,9)(-0.15,-0.11){16}{\line(-1,0){0.15}}
\multiput(67.82,7.20)(-0.24,-0.12){11}{\line(-1,0){0.24}}
\multiput(65.13,5.90)(-0.41,-0.11){7}{\line(-1,0){0.41}}
\multiput(62.24,5.17)(-1.49,-0.07){2}{\line(-1,0){1.49}}
\multiput(59.25,5.02)(-0.74,0.11){4}{\line(-1,0){0.74}}
\multiput(56.30,5.46)(-0.31,0.11){9}{\line(-1,0){0.31}}
\multiput(53.49,6.49)(-0.20,0.12){13}{\line(-1,0){0.20}}
\multiput(50.94,8.04)(-0.13,0.12){17}{\line(-1,0){0.13}}
\multiput(48.75,10.07)(-0.12,0.16){15}{\line(0,1){0.16}}
\multiput(47.01,12.50)(-0.11,0.25){11}{\line(0,1){0.25}}
\multiput(45.78,15.22)(-0.11,0.49){6}{\line(0,1){0.49}}
\put(45.12,18.13){\line(0,1){2.99}}
\multiput(45.04,21.12)(0.10,0.59){5}{\line(0,1){0.59}}
\multiput(45.56,24.06)(0.11,0.28){10}{\line(0,1){0.28}}
\multiput(46.65,26.84)(0.12,0.18){14}{\line(0,1){0.18}}
\multiput(48.27,29.35)(0.12,0.12){18}{\line(0,1){0.12}}
\multiput(50.36,31.49)(0.16,0.11){15}{\line(1,0){0.16}}
\multiput(52.83,33.17)(0.28,0.12){10}{\line(1,0){0.28}}
\multiput(55.58,34.33)(0.74,0.11){6}{\line(1,0){0.74}}
\multiput(100,35)(0.99,-0.10){3}{\line(1,0){0.99}}
\multiput(102.97,34.70)(0.36,-0.11){8}{\line(1,0){0.36}}
\multiput(105.83,33.82)(0.22,-0.12){12}{\line(1,0){0.22}}
\multiput(108.45,32.39)(0.13,-0.11){17}{\line(1,0){0.13}}
\multiput(110.74,30.47)(0.12,-0.15){16}{\line(0,-1){0.15}}
\multiput(112.60,28.14)(0.11,-0.22){12}{\line(0,-1){0.22}}
\multiput(113.96,25.48)(0.12,-0.41){7}{\line(0,-1){0.41}}
\multiput(114.77,22.60)(0.11,-1.49){2}{\line(0,-1){1.49}}
\multiput(115,19.63)(-0.09,-0.74){4}{\line(0,-1){0.74}}
\multiput(114.62,16.66)(-0.12,-0.35){8}{\line(0,-1){0.35}}
\multiput(113.67,13.83)(-0.11,-0.20){13}{\line(0,-1){0.20}}
\multiput(112.18,11.24)(-0.12,-0.13){17}{\line(0,-1){0.13}}
\multiput(110.20,9)(-0.15,-0.11){16}{\line(-1,0){0.15}}
\multiput(107.82,7.20)(-0.24,-0.12){11}{\line(-1,0){0.24}}
\multiput(105.13,5.90)(-0.41,-0.11){7}{\line(-1,0){0.41}}
\multiput(102.24,5.17)(-1.49,-0.07){2}{\line(-1,0){1.49}}
\multiput(99.25,5.02)(-0.74,0.11){4}{\line(-1,0){0.74}}
\multiput(96.30,5.46)(-0.31,0.11){9}{\line(-1,0){0.31}}
\multiput(93.49,6.49)(-0.20,0.12){13}{\line(-1,0){0.20}}
\multiput(90.94,8.04)(-0.13,0.12){17}{\line(-1,0){0.13}}
\multiput(88.75,10.07)(-0.12,0.16){15}{\line(0,1){0.16}}
\multiput(87.01,12.50)(-0.11,0.25){11}{\line(0,1){0.25}}
\multiput(85.78,15.22)(-0.11,0.49){6}{\line(0,1){0.49}}
\put(85.12,18.13){\line(0,1){2.99}}
\multiput(85.04,21.12)(0.10,0.59){5}{\line(0,1){0.59}}
\multiput(85.56,24.06)(0.11,0.28){10}{\line(0,1){0.28}}
\multiput(86.65,26.84)(0.12,0.18){14}{\line(0,1){0.18}}
\multiput(88.27,29.35)(0.12,0.12){18}{\line(0,1){0.12}}
\multiput(90.36,31.49)(0.16,0.11){15}{\line(1,0){0.16}}
\multiput(92.83,33.17)(0.28,0.12){10}{\line(1,0){0.28}}
\multiput(95.58,34.33)(0.74,0.11){6}{\line(1,0){0.74}}
\multiput(140,35)(0.99,-0.10){3}{\line(1,0){0.99}}
\multiput(142.97,34.70)(0.36,-0.11){8}{\line(1,0){0.36}}
\multiput(145.83,33.82)(0.22,-0.12){12}{\line(1,0){0.22}}
\multiput(148.45,32.39)(0.13,-0.11){17}{\line(1,0){0.13}}
\multiput(150.74,30.47)(0.12,-0.15){16}{\line(0,-1){0.15}}
\multiput(152.60,28.14)(0.11,-0.22){12}{\line(0,-1){0.22}}
\multiput(153.96,25.48)(0.12,-0.41){7}{\line(0,-1){0.41}}
\multiput(154.77,22.60)(0.11,-1.49){2}{\line(0,-1){1.49}}
\multiput(155,19.63)(-0.09,-0.74){4}{\line(0,-1){0.74}}
\multiput(154.62,16.66)(-0.12,-0.35){8}{\line(0,-1){0.35}}
\multiput(153.67,13.83)(-0.11,-0.20){13}{\line(0,-1){0.20}}
\multiput(152.18,11.24)(-0.12,-0.13){17}{\line(0,-1){0.13}}
\multiput(150.20,9)(-0.15,-0.11){16}{\line(-1,0){0.15}}
\multiput(147.82,7.20)(-0.24,-0.12){11}{\line(-1,0){0.24}}
\multiput(145.13,5.90)(-0.41,-0.11){7}{\line(-1,0){0.41}}
\multiput(142.24,5.17)(-1.49,-0.07){2}{\line(-1,0){1.49}}
\multiput(139.25,5.02)(-0.74,0.11){4}{\line(-1,0){0.74}}
\multiput(136.30,5.46)(-0.31,0.11){9}{\line(-1,0){0.31}}
\multiput(133.49,6.49)(-0.20,0.12){13}{\line(-1,0){0.20}}
\multiput(130.94,8.04)(-0.13,0.12){17}{\line(-1,0){0.13}}
\multiput(128.75,10.07)(-0.12,0.16){15}{\line(0,1){0.16}}
\multiput(127.01,12.50)(-0.11,0.25){11}{\line(0,1){0.25}}
\multiput(125.78,15.22)(-0.11,0.49){6}{\line(0,1){0.49}}
\put(125.12,18.13){\line(0,1){2.99}}
\multiput(125.04,21.12)(0.10,0.59){5}{\line(0,1){0.59}}
\multiput(125.56,24.06)(0.11,0.28){10}{\line(0,1){0.28}}
\multiput(126.65,26.84)(0.12,0.18){14}{\line(0,1){0.18}}
\multiput(128.27,29.35)(0.12,0.12){18}{\line(0,1){0.12}}
\multiput(130.36,31.49)(0.16,0.11){15}{\line(1,0){0.16}}
\multiput(132.83,33.17)(0.28,0.12){10}{\line(1,0){0.28}}
\multiput(135.58,34.33)(0.74,0.11){6}{\line(1,0){0.74}}
\put(8.33,10.33){\line(1,1){21.33}}
\put(6.67,13){\line(1,0){26.67}}
\put(47,13){\line(1,1){20.33}}
\put(48,10.67){\line(1,0){23.67}}
\put(113.67,13){\line(-1,0){27.33}}
\put(152,10.67){\line(-1,0){23.67}}
\put(107,33.33){\line(6,-5){4.67}}
\put(147,33.33){\line(6,-5){4.67}}
\put(39.67,20){\makebox(0,0)[cc]{$-$}}
\put(80.33,20){\makebox(0,0)[cc]{$=$}}
\put(120.67,20){\makebox(0,0)[cc]{$-$}}
\put(110,31.33){\line(0,-1){22.67}}
\put(152.67,12.33){\line(-1,6){3.22}}

\linethickness{1.5pt}

\multiput(10,9)(40,0){4}{\line(-2,3){3.78}}
\multiput(9.67,9)(40,0){4}{\line(-2,3){3.78}}

\multiput(33.67,13.67)(40,0){4}{\line(-3,-4){4}}
\multiput(33.33,13.67)(40,0){4}{\line(-3,-4){4}}

\multiput(26.67,33.33)(40,0){4}{\line(6,-5){5.67}}
\multiput(26.33,33.33)(40,0){4}{\line(6,-5){5.67}}

\thinlines

\end{picture}
\end{center}
\caption{The Four--Term relation}
\label{4T}
\end{figure}
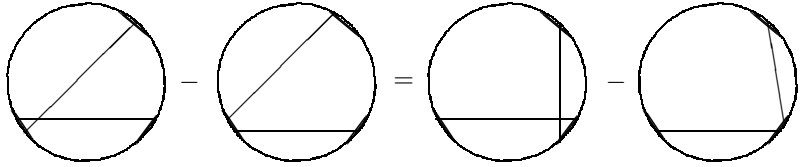

Consider the linear space of all chord diagrams
of degree $n$, where $n$ is a non--negative integer
over ${\bf Q}$.
Let us define the four--term (4T) relation as
follows, see Fig. \ref{4T}.

The summands represented in Fig. 1. are chord diagrams
with $n$ chords such that $n-2$ chords of all diagrams
have exactly the same position. The other two chords
are shown in the figure. They connect points lying
in the marked segments. There are no other vertices
on the marked segments except for those shown in the
figure.

Now, let us define the chord diagram algebra:

\begin{dfn}
The {\em chord diagram algebra} ${\cal A}^c$ is
a formal algebra, whose elements are equivalence
classes of linear combinations of chord diagrams
(possibly, having different degrees) modulo the
four--term relation.

The unity of this algebra is the equivalence
class of the chord diagram without chords.

The multiplication of two diagrams
$A$ and $B$ is defined as follows:
Let us break the circles of these diagrams at arbitrary
points which are not vertices,
and then connect them together according
to the circle orientation. An example of
two chord diagrams with marked points (left hand)
and their product (right hand) is shown in
Figure \ref{Prod}.

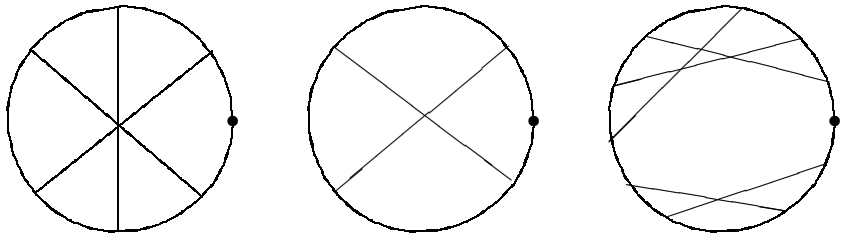
\begin{figure}
\unitlength 1mm
\linethickness{0.4pt}
\begin{center}
\begin{picture}(115.67,35)
\multiput(20,35)(0.99,-0.10){3}{\line(1,0){0.99}}
\multiput(22.97,34.70)(0.36,-0.11){8}{\line(1,0){0.36}}
\multiput(25.83,33.82)(0.22,-0.12){12}{\line(1,0){0.22}}
\multiput(28.45,32.39)(0.13,-0.11){17}{\line(1,0){0.13}}
\multiput(30.74,30.47)(0.12,-0.15){16}{\line(0,-1){0.15}}
\multiput(32.60,28.14)(0.11,-0.22){12}{\line(0,-1){0.22}}
\multiput(33.96,25.48)(0.12,-0.41){7}{\line(0,-1){0.41}}
\multiput(34.77,22.60)(0.11,-1.49){2}{\line(0,-1){1.49}}
\multiput(35,19.63)(-0.09,-0.74){4}{\line(0,-1){0.74}}
\multiput(34.62,16.66)(-0.12,-0.35){8}{\line(0,-1){0.35}}
\multiput(33.67,13.83)(-0.11,-0.20){13}{\line(0,-1){0.20}}
\multiput(32.18,11.24)(-0.12,-0.13){17}{\line(0,-1){0.13}}
\multiput(30.20,9)(-0.15,-0.11){16}{\line(-1,0){0.15}}
\multiput(27.82,7.20)(-0.24,-0.12){11}{\line(-1,0){0.24}}
\multiput(25.13,5.90)(-0.41,-0.11){7}{\line(-1,0){0.41}}
\multiput(22.24,5.17)(-1.49,-0.07){2}{\line(-1,0){1.49}}
\multiput(19.25,5.02)(-0.74,0.11){4}{\line(-1,0){0.74}}
\multiput(16.30,5.46)(-0.31,0.11){9}{\line(-1,0){0.31}}
\multiput(13.49,6.49)(-0.20,0.12){13}{\line(-1,0){0.20}}
\multiput(10.94,8.04)(-0.13,0.12){17}{\line(-1,0){0.13}}
\multiput(8.75,10.07)(-0.12,0.16){15}{\line(0,1){0.16}}
\multiput(7.01,12.50)(-0.11,0.25){11}{\line(0,1){0.25}}
\multiput(5.78,15.22)(-0.11,0.49){6}{\line(0,1){0.49}}
\put(5.12,18.13){\line(0,1){2.99}}
\multiput(5.04,21.12)(0.10,0.59){5}{\line(0,1){0.59}}
\multiput(5.56,24.06)(0.11,0.28){10}{\line(0,1){0.28}}
\multiput(6.65,26.84)(0.12,0.18){14}{\line(0,1){0.18}}
\multiput(8.27,29.35)(0.12,0.12){18}{\line(0,1){0.12}}
\multiput(10.36,31.49)(0.16,0.11){15}{\line(1,0){0.16}}
\multiput(12.83,33.17)(0.28,0.12){10}{\line(1,0){0.28}}
\multiput(15.58,34.33)(0.74,0.11){6}{\line(1,0){0.74}}
\multiput(60,35)(0.99,-0.10){3}{\line(1,0){0.99}}
\multiput(62.97,34.70)(0.36,-0.11){8}{\line(1,0){0.36}}
\multiput(65.83,33.82)(0.22,-0.12){12}{\line(1,0){0.22}}
\multiput(68.45,32.39)(0.13,-0.11){17}{\line(1,0){0.13}}
\multiput(70.74,30.47)(0.12,-0.15){16}{\line(0,-1){0.15}}
\multiput(72.60,28.14)(0.11,-0.22){12}{\line(0,-1){0.22}}
\multiput(73.96,25.48)(0.12,-0.41){7}{\line(0,-1){0.41}}
\multiput(74.77,22.60)(0.11,-1.49){2}{\line(0,-1){1.49}}
\multiput(75,19.63)(-0.09,-0.74){4}{\line(0,-1){0.74}}
\multiput(74.62,16.66)(-0.12,-0.35){8}{\line(0,-1){0.35}}
\multiput(73.67,13.83)(-0.11,-0.20){13}{\line(0,-1){0.20}}
\multiput(72.18,11.24)(-0.12,-0.13){17}{\line(0,-1){0.13}}
\multiput(70.20,9)(-0.15,-0.11){16}{\line(-1,0){0.15}}
\multiput(67.82,7.20)(-0.24,-0.12){11}{\line(-1,0){0.24}}
\multiput(65.13,5.90)(-0.41,-0.11){7}{\line(-1,0){0.41}}
\multiput(62.24,5.17)(-1.49,-0.07){2}{\line(-1,0){1.49}}
\multiput(59.25,5.02)(-0.74,0.11){4}{\line(-1,0){0.74}}
\multiput(56.30,5.46)(-0.31,0.11){9}{\line(-1,0){0.31}}
\multiput(53.49,6.49)(-0.20,0.12){13}{\line(-1,0){0.20}}
\multiput(50.94,8.04)(-0.13,0.12){17}{\line(-1,0){0.13}}
\multiput(48.75,10.07)(-0.12,0.16){15}{\line(0,1){0.16}}
\multiput(47.01,12.50)(-0.11,0.25){11}{\line(0,1){0.25}}
\multiput(45.78,15.22)(-0.11,0.49){6}{\line(0,1){0.49}}
\put(45.12,18.13){\line(0,1){2.99}}
\multiput(45.04,21.12)(0.10,0.59){5}{\line(0,1){0.59}}
\multiput(45.56,24.06)(0.11,0.28){10}{\line(0,1){0.28}}
\multiput(46.65,26.84)(0.12,0.18){14}{\line(0,1){0.18}}
\multiput(48.27,29.35)(0.12,0.12){18}{\line(0,1){0.12}}
\multiput(50.36,31.49)(0.16,0.11){15}{\line(1,0){0.16}}
\multiput(52.83,33.17)(0.28,0.12){10}{\line(1,0){0.28}}
\multiput(55.58,34.33)(0.74,0.11){6}{\line(1,0){0.74}}
\multiput(100,35)(0.99,-0.10){3}{\line(1,0){0.99}}
\multiput(102.97,34.70)(0.36,-0.11){8}{\line(1,0){0.36}}
\multiput(105.83,33.82)(0.22,-0.12){12}{\line(1,0){0.22}}
\multiput(108.45,32.39)(0.13,-0.11){17}{\line(1,0){0.13}}
\multiput(110.74,30.47)(0.12,-0.15){16}{\line(0,-1){0.15}}
\multiput(112.60,28.14)(0.11,-0.22){12}{\line(0,-1){0.22}}
\multiput(113.96,25.48)(0.12,-0.41){7}{\line(0,-1){0.41}}
\multiput(114.77,22.60)(0.11,-1.49){2}{\line(0,-1){1.49}}
\multiput(115,19.63)(-0.09,-0.74){4}{\line(0,-1){0.74}}
\multiput(114.62,16.66)(-0.12,-0.35){8}{\line(0,-1){0.35}}
\multiput(113.67,13.83)(-0.11,-0.20){13}{\line(0,-1){0.20}}
\multiput(112.18,11.24)(-0.12,-0.13){17}{\line(0,-1){0.13}}
\multiput(110.20,9)(-0.15,-0.11){16}{\line(-1,0){0.15}}
\multiput(107.82,7.20)(-0.24,-0.12){11}{\line(-1,0){0.24}}
\multiput(105.13,5.90)(-0.41,-0.11){7}{\line(-1,0){0.41}}
\multiput(102.24,5.17)(-1.49,-0.07){2}{\line(-1,0){1.49}}
\multiput(99.25,5.02)(-0.74,0.11){4}{\line(-1,0){0.74}}
\multiput(96.30,5.46)(-0.31,0.11){9}{\line(-1,0){0.31}}
\multiput(93.49,6.49)(-0.20,0.12){13}{\line(-1,0){0.20}}
\multiput(90.94,8.04)(-0.13,0.12){17}{\line(-1,0){0.13}}
\multiput(88.75,10.07)(-0.12,0.16){15}{\line(0,1){0.16}}
\multiput(87.01,12.50)(-0.11,0.25){11}{\line(0,1){0.25}}
\multiput(85.78,15.22)(-0.11,0.49){6}{\line(0,1){0.49}}
\put(85.12,18.13){\line(0,1){2.99}}
\multiput(85.04,21.12)(0.10,0.59){5}{\line(0,1){0.59}}
\multiput(85.56,24.06)(0.11,0.28){10}{\line(0,1){0.28}}
\multiput(86.65,26.84)(0.12,0.18){14}{\line(0,1){0.18}}
\multiput(88.27,29.35)(0.12,0.12){18}{\line(0,1){0.12}}
\multiput(90.36,31.49)(0.16,0.11){15}{\line(1,0){0.16}}
\multiput(92.83,33.17)(0.28,0.12){10}{\line(1,0){0.28}}
\multiput(95.58,34.33)(0.74,0.11){6}{\line(1,0){0.74}}
\put(35,19.67){\circle*{1.33}}
\put(75,19.67){\circle*{1.33}}
\put(115,19.67){\circle*{1.33}}
\multiput(32.33,29)(-0.15,-0.12){156}{\line(-1,0){0.15}}
\put(19.67,35){\line(0,-1){30}}
\multiput(31,9.67)(-0.14,0.12){162}{\line(-1,0){0.14}}
\put(71.67,29.67){\line(-6,-5){23}}
\put(48.33,29.67){\line(4,-3){23.67}}
\put(114.33,25){\line(-4,1){24.33}}
\put(110.67,30.67){\line(-4,-1){25}}
\put(102.67,34.67){\line(-1,-1){17.67}}
\put(87.33,11.33){\line(6,-1){21.33}}
\put(92.67,7){\line(3,1){21}}

\end{picture}
\end{center}

\caption{Product of Chord Diagrams}
\label{Prod}
\end{figure}

Thus we obtain the chord diagram that can
be treated as the product $AB$. However,
this product is not uniquely defined;
it depends on the choice of break points.
Bar--Natan \cite{BN} showed that this
choice of fixed point does not change
the final result up to the 4T--relation.

Obviously, the product of two linear combinations
of chord diagrams is defined according to distributivity
rule

\end{dfn}

Let us denote by
${\cal A}^c_{n}$ the linear space of chord
diagrams of degree $n$ modulo $4T$--relation.

Now, let us define another algebra

\begin{dfn}
A {\em Chinese character diagram} is a cubic graph, with all vertices either exterior or interior,
with indicated oriented cycle, containing all exterior vertices. At each interior vertex
one also indicates a cyclic order of outcoming edges counterclockwise.
Chinese character diagrams are considered
up to graph equivalece preserving the
oriented cycle ({\em circle}).

The degree of a CCD is half the number of
its vertices
\end{dfn}

\begin{re}

In the sequel, all circles shown
in figures are thougth to be oriented
counterlockwise; all structures at
interior points are taken from the
plane.

Although the degree of a CCD is defined by its number of vertices divided by two, it turns out that it is an integer, as can easily been shown.

\end{re}

Note, that each chord diagram is
a Chinese character diagram without
exterior vertices. In this case,
its CD degree coincides with its
CCD degree.

Chinese character diagrams admit a similar
algebraic structure.

\begin{figure}
\unitlength 0.4mm
\begin{center}
\begin{picture}(100,60)

\thicklines
\put(40,10){\line(-1,0){30}}
\put(90,10){\line(-1,0){30}}
\put(140,10){\line(-1,0){30}}

\thinlines
\put(25,10){\line(0,1){20}}
\put(25,30){\line(1,1){20}}
\put(25,30){\line(-1,1){20}}

\put(120,10){\line(2,3){25}}
\put(130,10){\line(-2,3){25}}

\put(80,10){\line(2,5){15}}
\put(70,10){\line(-2,5){15}}

\put(45,30){\Large{=}}
\put(100,30){\Large{-}}

\end{picture}
\end{center}
\caption{ $STU$--relation}
\label{STU}

\end{figure}
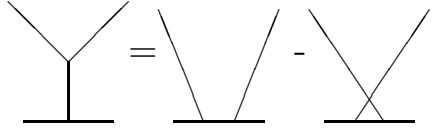

Consider the space of linear combinations
of Chinese character diagrams. Let us factorize it
by the so--called {\em
$STU$--relation}, see Fig. \ref{STU}.
Denote the obtained space by
${\cal A}^t$.

As before, we have shown only a small changing
part of the CCD. All other vertices remain the
same in all three cases shown above.

By using the $STU$--relation, each CCD
can be transformed into a linear combination of
chord diagrams.

Indeed, one can resolve an interior vertex connected
with an exterior vertex. Thus, if the number of
interior vertices is greater than 0, it can be
decreased. The final result (i.e. a linear combination
of chord diagrams) is uniquely defined up to
$4T$--relation, see \cite{BN}.
Thus, the space ${\cal A}^t$ inherits
the algebraic structure of
${\cal A}^c$.

Let us take into account the two relations
that take place in ${\cal A}^{t}$, see \cite{BN}.

\begin{th}
Following identities hold in ${\cal A}^{t}$:
\begin{enumerate}

\item Antisymmetry at interior points, see Fig.\ref{Anti}

\begin{figure}
\unitlength 0.6mm
\begin{center}
\begin{picture}(100,60)
\thinlines
\put(30,10){\line(0,1){20}}
\put(30,30){\line(1,1){20}}
\put(30,30){\line(-1,1){20}}
\put(80,10){\line(0,1){20}}
\put(80,40){\oval(20,20)[b]}
\put(70,40){\line(3,1){30}}
\put(90,40){\line(-3,1){30}}

\put(50,30){\Large{+}}
\put(100,30){\Large{=0}}

\end{picture}
\end{center}
\caption{The antisymmetry relation}
\label{Anti}

\end{figure}
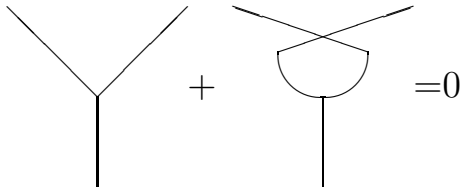

\item The $IHX$--relation, see Fig.\ref{IHX}

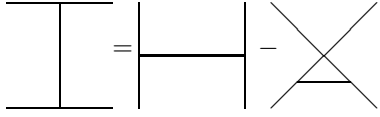
\begin{figure}
\begin{center}
\begin{picture}(100,60)
\thinlines
\put(10,10){\line(1,0){40}}
\put(30,10){\line(0,1){40}}
\put(10,50){\line(1,0){40}}

\put(60,10){\line(0,1){40}}
\put(100,10){\line(0,1){40}}
\put(60,30){\line(1,0){40}}

\put(110,10){\line(1,1){40}}
\put(150,10){\line(-1,1){40}}
\put(120,20){\line(1,0){20}}

\put(50,30){$\Large{=}$}
\put(105,30){$\Large{-}$}

\end{picture}
\end{center}
\caption{IHX--relation}
\label{IHX}
\end{figure}

\end{enumerate}

\end{th}

Both ${\cal A}^{c}$ and ${\cal A}^{t}$ are
graded according to the half number of the
total quantity of vertices
(for chord diagrams the degree equals the number
of chords);
the factorization described above preserves
this structures because
both  the $4T$-- and $STU$--
relations are homogeneous according to this graduation.

\begin{figure}
\unitlength 0.8mm
\begin{center}
\begin{picture}(75,35)
\thicklines
\put(41.33,15){\oval(67.33,20)[b]}
\thinlines
\put(25,5){\line(0,1){10}}
\put(35,5){\line(0,1){10}}
\put(45,5){\line(0,1){10}}
\put(55,5){\line(0,1){10}}
\put(23.33,17.67){\makebox(0,0)[cc]{$x_1$}}
\put(33.67,17.67){\makebox(0,0)[cc]{$x_2$}}
\put(44,17.67){\makebox(0,0)[cc]{$x_3$}}
\put(53.67,17.67){\makebox(0,0)[cc]{$x_4$}}
\put(19,2){\makebox(0,0)[cc]{$u$}}
\put(85.33,2.33){\makebox(0,0)[cc]{$R(x_4)R(x_3)R(x_2)R(x_1)u$}}
\thicklines
\put(41.17,16.67){\oval(67,33.33)[t]}
\end{picture}
\end{center}
\caption{Contraction}
\label{Contr}
\end{figure}
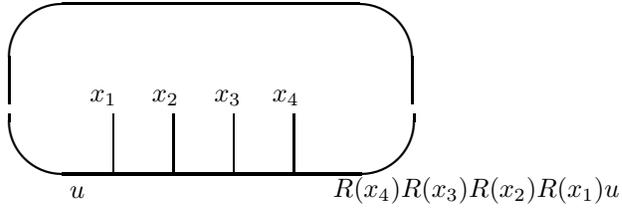

As it is  shown in \cite{BN}, by using a given representation of a
semisimple Lie algebra $G$ and a $n$--Chinese character diagram,
one can obtain a tensor of type $(n,0)$ invariant under the adjoint action
of $G$. This can be done as follows: Let $R$ be a representation of some Lie algebra $G$. Consider the $n$--Chinese diagram $C$
without interior vertices (i.e. all trivalent vertices are exterior).
One defines the $n$--linear form on "tails", corresponding to $C$ by
$f(x_{1},\dots,x_{n})=Tr (R(x_{1})\dots R(x_{n}))$,
see figure \ref{Contr}; where $TR$ is the trace of an operator acting on the Lie algebra.

The computation of the multilinear form for arbitrary $n$--chinese
diagram can be done by contracting "tails" at interior vertices by using the trivalent structural
tensor $c^{i}_{jk}$ and the bivalent metrics $g_{ij}$.
Fix a Lie algebra $G$ and its representation $R$. The map from the set of  $k$--Chinese character
diagrams to the adjoint invariant $k$--linear forms on $G$ is invariant under the $STU$--relation. If
a linear combination of $k$--Chinese character diagrams equals zero, up to the $STU$--relation,
then the corresponding tensor equals zero by the Jacobi identity and the skew symmetry, see \cite{BN}.
According to the isomorphism, ${\cal A}^t$ ¨ ${\cal A}^c$ all $k$--linear forms, coming from $k$--Chinese
character diagrams, can be obtained only from $k$--chord diagrams.
Later we shall give an algorithm for calculating these tensors for some representations of Lie algebras
and arbitrary $k$--chord diagrams.

Note that the case of chord diagrams without tails is also interesting. In this case we get scalars, which are
invariant under the 4T--relation. These invariants give some properties of the chord diagram algebra. As shown
in \cite{BN}, this algebra plays a significant role in the Vassiliev invariant theory.

\section{Adjoint representation of sl(n) and so(n)}

\subsection{The case of $sl(n)$}

Consider the $k$--Chinese character diagram $L(k)$, consisting of a circle and $k$ outcoming
tails. In the sequel, we shall denote it simply by $L$. As told in the introduction, for such a diagram
we get the following invariant tensor:

$$f(x_{1},\dots,x_{k})=Tr (ad_{x_{1}}\dots ad_{x_{k}})\eqno(1)$$

So, we have to find the operator trace $O=ad_{x_1}\dots ad_{x_{k}}$ for $sl(n)$.
First, consider this trace for $gl(n)$.

\begin{lm}
The final result after a contraction is the same for both $sl(n)$ and $gl(n)$.
\end{lm}

\begin{proof}
Any contraction means a calculation of a trace for an operator, acting
on $n\times n$ matrices. For these matrices, choose the basis
consisting of the unit matrix and of matrices of $sl(n)$. The initial operator $O$
is a composition of commutation operators, so it vanishes on the identitity matrix. Thus,
all operators obtained from it vanish on the identity as well. Thus, for each
of them, the traces in bases of $sl(n)$ and $gl(n)$ coincide.
\end{proof}

Choose the basis $E_{ij}, i,j=1,\dots, n$ of $gl(n)$, where $E_{ij}$ is the matrix having a $1$ in the $i^{th}$ row and $j^{th}$ columns and zero elsewhere. The dual basis is $E_{ji},i,j=1,\dots, n$. The formula
(1) looks like

$$f(x_{1},\dots,x_{k})=\sum_{i,j=1}^{n} Tr (E_{ij} ad_{x_1}\dots ad_{x_n} E_{ij})\eqno(2),$$
where $Tr$ is just the matrix trace.

In the sequel, for the case of the adjoint representation, $Tr$ is the usual matrix trace.

Now rewrite (2) representing
the action $ad_{x_i}$ on $u$ as a commutator. So,

$$f(x_{1},\dots, x_{k})=
\sum_{{l=0} \atop {i_{1}>\dots >i_{l},i_{l+1}<\dots<i_{n}}} ^{n}
(-1)^{l}
Tr (x_{i_1}\dots x_{i_{k}}) Tr (x_{i_{k+1}}\dots x_{i_n}). \eqno(3). $$

That is, indices of $x$
before $u$ are descending and
those after $u$ are ascending.
The sum is taken according to all the
decompositions of indices $1,\dots,k$
into "left ones" and "right ones"

\begin{re}
For the empty diagram ($k=0$) the formula (3) does not hold.
The corresponding scalar equals $n^{2}-1$.
\end{re}

The two following types of members must be determined:

$$ \sum_{\alpha} Tr (P \alpha Q \tilde\alpha)\eqno(4)$$
and
$$ \sum_{ij} Tr (P \alpha) Tr (Q \tilde\alpha)\eqno(5),$$
where $P,Q \in gl(n)$,
and $\alpha, \tilde\alpha$ run the usual and dual bases of the Lie algebra. Therefore equations (4) and (5) take the form:

$$ \sum_{ij} Tr (P E_{ij} Q E_{ji})= \sum_{i,j=1}^{n} P_{jj} Q_{ii}= Tr(P) Tr(Q) \eqno(6)$$
$$ \sum_{ij} Tr (P E_{ij}) Tr (Q E_{ji}) = \sum_{i,j=1}^{n} P_{ji}Q_{ij}= Tr (P Q) \eqno(7).$$

Computing these contractions directly for
$sl(n)$ would result in a much more
complicated algorithm.

\begin{lm}
For $sl(n)$ the sum (4) equals $$Tr(P) Tr(Q)- {1\over n} Tr (PQ),\eqno(*)$$
and (5) equals $$Tr (PQ)-  {1\over n} Tr(P) Tr(Q).\eqno(**)$$
\end{lm}

\begin{proof}

Choose the basis consisting of $E_{ij},i\neq j$
and coadjoint elements
$d_{l}=E_{11}+\dots E_{l-1,l-1}-(l-1)E_{ll}, l=1,\dots {n-1}$.
of $sl(n)$.
The dual basis consists of $E_{ji}$ and
$d_{l}=E_{11}+\dots E_{l-1,l-1}-(l-1)E_{ll}, l=1,\dots {n-1}$.

We have:

$$ \sum_{i \neq j} Tr (P E_{ij} Q E_{ji})=
\sum_{i,j=1}^{n} P_{jj} Q_{ii}-\sum_{i=1}^{n} P_{ii} Q_{ii}
= $$

$$Tr P Tr Q - \sum_{i=1}^{n} P_{ii} Q_{ii}$$

$$
\sum_{l=1}^{n-1} Tr (P d_{l} Q d_{l})\frac{1}{l^{2}-l}$$

$$=\left({1\over 2}+{1\over 6}+{1\over 12}+\dots +
+{1\over {n(n-1)} }\right)\sum_{i=1}^{n} P_{ii}Q_{ii}+$$

$$\left({-{1\over 2}+{1\over 6}+{1\over 12}+\dots+ {1\over n(n-1)}}\right)
\sum_{i\neq j}P_{ij}Q_{ji}=-{1\over n}Tr (P Q)+ \sum_{i=1}^{n} P_{ii}Q_{ii}$$

Collecting these sums together, we get:
$$Tr(P) Tr(Q) - {1\over n} Tr (P Q)$$.

$$ \sum_{i \neq j} Tr (P E_{ij}) Tr Q (E_{ji}) =
\sum_{i,j=1}^{n} P_{ji}Q_{ij}= Tr (P Q) - \sum_{i=1}^{n} P_{ii}Q_{ii} .$$

$$\sum_{i=1}^{n-1} Tr (P d_{l}) Tr (P d_{l}) \frac{1}{l^{2}-l} =$$

$$\left(-{1\over 2}+ {1\over 6}+ {1\over 12}+ \dots {1\over {n(n-1)}}\right)
\sum_{i\neq j}P_{ii}Q_{jj}+\sum_{i=1}^{n} P_{ii}Q_{ii}$$

Finally, we get:
$$Tr(P) Tr(Q) -{1\over n} Tr (PQ)$$.

\end{proof}

So, in order to calculate the tensor
corresponding to a chord diagram of order $k$
for the case of the adjoint representation
of $so(n)$, we have to calculate
$2^{2k}$ summands. Since we consider
the trace in $gl(n)$, the number
of summands increases while contracting.
We get the following

\begin{th}
Following conditions hold:
\begin{enumerate}
\item Each $m$--linear form on $x_{i_1},\dots x_{i_m}$, obtained
from the coadjoint representation of $sl(n)$
by contracting the $k$--linear form corresponding to
$L$, is a linear combination of traces of
products for $x_{i_j}$ and $x_{i_j}^{*}$ with
polynomial (with respect to $n$) coefficients.

\item the power of these polynomials does not exceed $k+2$.
\end{enumerate}
\end{th}

Note, that the statement of the theorem
for Chinese character diagrams without
univalent vertices was proved by
Bar--Natan, \cite{BN2}. Here we
give our own

\begin{proof}
Consider the desired tensor as a sum of
$2^{2k}$ summands according to (3) and contract it in $gl(n)$.
For each of $2^{2k}$ summands we use the induction
on the number $q$ of contracted elements.
Let us also observe the number of multipliers (traces) in
this product.
For $q=0$ the statement is evident.
Consider the power of a summand and the number of its
multipliers.
In the case of a contraction like (6), one gets either

1. {$P=Q=E$}. The power increases by two, the number of multipliers decreases by 1.

2. {Either $P$ or $Q$ equals the identity matrix}.
The power decreases by 1, the number of multipliers remains the same.

3. {None of $P$ and $Q$ is the identity matrices}.
The power stays the same and the number of multipliers increases by 1.

In the case of the contraction (7), one of the following
situation occurs:

1.{$P=Q=E$}. The power stays the same, the number of multipliers decreases by 1.

2.{At least one of $P,Q$ is not equal to $E$}.
The power stays the same, the number of multipliers decreases by 1.

Thus we see that contraction implies the induction step
for the statement 1). Besides, after $l$ contractions the
sum of coefficient's maximal power for a summand and the number
of multipliers of this summand does not exceed $l+2$.
Since the number of multipliers is not negative and the number of
contractions satisfies $l\le k$, we obtain that the power of a coefficient for a
summand does not exceed $k+2$.
\end{proof}

\subsection{The case of $so(n)$}

As in the case of $sl(n)$, we begin with the $k$--Chinese
diagram $L$ consisting of a circle with $k$ outcoming tails.
To find (3), we have to calculate (4) and (5).

\begin{lm}
For $so(n)$ the sum (4) looks like
$${1\over 2} (Tr(P) Tr(Q)- Tr (PQ^{*}))\eqno(8),$$
and (5) looks like $${1\over 2}(Tr (PQ)- Tr (P) Tr (Q^{*}))\eqno(9)$$
\end{lm}

\begin{proof}
Choose the selfadjoint basis $E_{ij}-E_{ji},i> j$
of so(n), whose elements have length $2$. So:

$$\sum_{i>j} {1\over 2} Tr (P (E_{ij}-E_{ji}) Q (E_{ji}-E_{ij})=$$

$${1\over 2} \sum_{i>j} Tr P [ P_{jj}E_{ii}+ P_{ii}E_{jj}- Q_{ji} E_{ij}- Q_{ij} E_{ji}]= $$

$$={1\over 2}\sum_{i>j} (P_{ii} Q_{jj}- P_{ji} Q_{ji}- P_{ij}Q_{ij}=$$

$$ {1\over 2} (Tr (P) Tr (Q) - \sum_{i>j} (P_{ij}Q_{ij}+P_{ji}Q_{ji})+\sum_{i} P_{ii}Q_{ii}))=$$

$$= {1\over 2} (Tr (P) Tr(Q) - Tr (PQ^{*})).$$

$$\sum_{i>j} Tr (P(E_{ij}-E_{ji})) Tr (Q (E_{ji}-E_{ij}))=$$

$$\sum_{i>j} {1\over 2}((P_{ji}Q_{ij}+P_{ij}Q_{ji})-(P_{ii}Q_{ij}+P_{ji}Q_{ji}))=$$

$${1\over 2} (Tr(PQ)-Tr(PQ^{*})).$$
\end{proof}

From the formulae above we rewrite equation (3) as:
\begin{eqnarray*}
f(x_{1},\dots, x_{k})=
{1\over 2}\sum_{{l=0} \atop {i_{1}>\dots >i_{l},i_{l+1}<\dots<i_{n}}} ^{n}
(-1)^{l}\times
(Tr (x_{i_1}\dots x_{i_{k}}) Tr (x_{i_{k+1}}\dots x_{i_n})+\\
Tr(x_{i_1}\dots x_{i_{k}} x_{i_n}^{*}\dots x_{i_{k+1}}^{*}))
\end{eqnarray*}

\begin{th}
Each $m$--linear form on $x_{i_1},\dots x_{i_m}$, obtained from
the adjoint representation for $so(n)$ is a
linear combination of traces of products for $x_{i_j}$ and $x_{i_j}^{*}$
with polynomial coefficients with respect to $n$.
\end{th}

The proof is essentially the same as for $sl(n)$ with
formulae (8), (9) instead of (6) and (7).

\section{Arbitrary irreducible representations of sl(n)}

Note that the calculations for arbitrary representations of $sl(2)$ were done
in \cite{CV}.

In the case of arbitrary representation $R$ of ܬ $sl(n)$, it is important to caculate the formulae for
$\sum Tr (P\alpha Q\tilde \alpha)$  and  $\sum Tr (P \alpha) T(Q \tilde \alpha)$,
where $P$ and $Q$ are products of certain number of matrices,
representing $x_{i}$, $\alpha=R(\omega)$ and $\tilde \alpha= R(\tilde\omega)$,
where $\omega$ runs a basis of $sl(n)$ and $\tilde \omega$ runs the dual basis.

It is obvious that $\sum Tr (P\alpha Q\tilde \alpha)=A Tr P Tr Q + B Tr (PQ),$
where $A$ and $B$ are numbers depending only on $n$.
One finds them by using the undefined coefficients method. For $P=Q=E$, one gets

$$Dim(R) Dim(G)= A Dim(R)^2+ B Dim(G) \mbox{ or } Dim(G)= A Dim(R)+ B \eqno(10)$$

The second equation for calculating $A,B$ one obtains assuming $P=Q=R(H),$
where $H$ is an element of the Cartan subalgebra of $G$.
The basis ${\omega}$ and the dual basis ${\tilde \omega}$
consist of vectors of the Cartan subalgebra and root vectors.

So,

$$ \sum _{\alpha} Tr (R(H)\alpha R(H)\tilde \alpha)=
\sum_{1}+\sum_{2}. $$

Here $\sum_{1}, \sum_{2}$ are sums with respect to
$\alpha$ that runs all simple roots and all root vectors.
$  \sum_{1}=  \sum Tr(R(h)^2 \alpha \tilde \alpha)  $,
$  \sum_{2}= \sum Tr(R(H)R(e_\beta)R(H) R(e_{-\beta})) $

Taking into account

$$H e_\beta=e_\beta H+ (\beta,H)e_\beta$$ and

$$H ¥_{-\beta} =e_{-\beta} H- (\beta,H) e_{-\beta},$$
where $\beta$ is an arbitrary root of $G$, we get:

$$Tr \sum_{2}= Tr \sum_{2} R(H) e_\beta R(H) e_{\beta} =
R(H)^2 \sum_{2} Tr R(H) R(H)_\beta= \sum \beta(H)^2.$$
The last sum is taken with respect to the set of positive roots  $\beta$ of $G$.

Thus we obtain $$\sum \alpha R(H) \tilde \alpha R(H)= Dim(R) \sum \alpha \tilde \alpha+
+\sum b(H)^2 = B Tr R(H)^2 \eqno(11)$$

Equations (10),(11) are sufficient to find $A,B$.

In the case of the irreducible $k+1$--dimensional representation $R$ of $G=sl(2)$, (10),(11) result in
$A=\frac {2}{k+1}, B=\frac{1-2k}{1+k} $.

\section{The case of the adjoint representation of sl(n) and invariants
of the chord diagram algebra}

Consider a chord diagram $C$ and construct the corresponding scalar
,depending on $n$, by using the adjoint representation of
$sl(n)$. By theorem 2, this scalar is a polynomial
with respect to $n$; its degree does not exceed $k+2$, where $k$ is
the order of $C$. Besides, this polynomial is an invariant under
$4T$--relation. Denote this function by $U(C)$.

Now, let us find those chord diagrams $C$ for which
the power of $U(C)$ with respect to $n$ equals
$k+2$. Consider the decomposition (3) for a tensor,
corresponding to $C$, and try to filter those
diagrams, generating after ${\pm}n^{k+2}$ after complete contraction.
Fix a summand $S$ and observe its contraction. As it is shown in the
proof of theorem 2, any contraction of type (6) increases the sum of
coefficient power and the number of multipliers (traces) by $1$,
and any  contraction of type (7) decreases this sum.
Taking into account that after all contractions the number of
multipliers equals zero we see that, in order to obtain the
maximal degree $k+2$, we always have to contract by (6).

This means that for each contraction chord ends should belong to the same
multiplier.

\begin{dfn}
A {\em $d$--diagram} is a chord diagram whose set of chords can be split into two families
of pairwise nonintersecting chords.
\end{dfn}

\begin{th}
The power of $U(C)$ equals $k+2$ on a diagram $C$ on $k$ chords iff $C$ is a $d$--diagram.
Moreover, the coefficient at $n^{k+2}$ of $U(C)$ equals the number of splitting chords of
$C$ into two families (the first and the second) of pairwise nonintersecting
chords.
\end{th}

\begin{proof}
Consider the diagram $C$. As it is shown above, to obtain $n^{k+2}$
by contracting some summand $S$, the ends of contracted chord of the diagram
$C$ should belong to the same multiplier.

This means that for each step any chord must belong to the same multiplier with both ends
of it, since otherwise the contraction (7) does not allow to obtain the maximal possible power.

Consider some splitting of our chord diagram vertices
into some nonintersecting subsets in such a way that
both ends of each chord lie in the same subset.
In this case, one can say that each chord belongs
to some subset.
If we contract along a chord lying in one subset
(contraction $Tr(P\alpha Q\tilde \alpha) \to Tr P Tr Q$
of type (6)), then the subset containing $P$ and $Q$
is decomposed into two subsets: $P$ and $Q$.
If we want each chord this splitting to have both ends
lying in the same family, each two chords of the same family
of the initial diagram must not intersect each other.
The inverse statement is also true: if any two chords of the
same family for the initial diagram do not intersect each other,
then, by contracting a chord by (6), we get a diagram of smaller
degree where each two chords of the same family do not intersect each
other.
In (3) we have $2^{2k}$ members; each of them
corresponds to a splitting of chords of the diagram
into two families. Chords in each family must be pairwise
nonintersecting, i.e. $C$ is a $d$--diagram.
To conclude the proof of the second statement, we only have to
consider the coefficient at $k+2$th degree of the polynomial for
our $d$--diagram $C$. Among $2^{2k}$ members of (3) we have some number
of "good" members, giving $n^{k+2}$ after contraction.
Their number equals the number of splittings of chords
into two families in a proper way. After the final contraction,
each of them gives $n^{k+2}\cdot (-1)^{l}$, where $l$ --- is
the number of chords in one of two families. Since both ends of each
chord lie in the same family, $(-1)^{l}=1$, that concludes the proof.
\end{proof}

\begin{th} For each natural $n$ in each basis
of the space ${\cal A}^{c}_n$ contains at least one $d$--diagram.
\end{th}

\begin{proof}
Consider a basis of ${\cal A}^{c}_{k}$ consisting of some diagrams
$B_{1},\dots, B_{r}$. Suppose that there is no $d$--diagram among them.
Then the degree of the polynomial $\forall i=1,\dots r \quad U(B_{i})$
is less than $k+2$. Consider any $d$--diagram $C$ on $k$ chords. Since
$B_{1},\dots, B_{r}$ form a basis of ${\cal A}^{c}_{k}$, then
$C$ must be represented by a linear combination
of $B_{1},\dots, B_{r}$ modulus the $4T$--relation. Since $U(\cdot)$ is invariant
under the $4T$--relation, $U(C)$ is a linear combination of $U(B_{i}), i=1,\dots r$,
that is impossible.
\end{proof}

\begin{th}
For a chord diagram $C_k$ on $k$ chords all monomials of $U(C_{k})$ have the same parity as $n$.
\end{th}

\begin{proof}
For each contraction of a monomial the sum of its coefficient's power
and the number of its multipliers either decrease by 1 or increase by 1.
Since this sum was first equal to 2, then, after $n$ contractions when the
number of multipliers equals zero, the parity of the monomial power equals $n$.
\end{proof}

\begin{th}
For an arbitrary chord diagram $C$
the polynomial $U(C)$ is divisible by $n^{2}-1$, where this number equals the dimension of $sl(n)$. Moreover, for arbitrary chord diagrams
$C_1$ and $C_2$ the following formula holds:

$$U({C_1}\cdot {C_2})= {{U({C_1})\cdot U({C_2})}\over (n^{2}-1)}$$
\end{th}

\begin{proof}
Consider the chord diagram $C$ of degree $k$
and choose a chord $a$ of it.
Consider the formula (3) for the diagram $L$.
For each summand of it, let us contract all chords
of it except for the chord $a$. Thus we obtain
a sum of monomials corresponding to
the $1$--chord diagram
obtained from $C$ by breaking the chord $a$.

Now, let us contract the rest in the sence of
$gl(n)$.

In the case
(*) we get the additional coefficient
$Tr(E)Tr(E)-{1\over n} Tr(E\cdot E)= n^{2}-1$.
In the case
(**) we get the coefficient $Tr(E\cdot E)-{1\over n} Tr(E\cdot E)=0$.

Thus, while performing the last step in the sense of
$gl(n)$, each coefficient either vanishes or is multiplied
by $(n^{2}-1)$. That completes the proof of the first part
of the theorem.

Let us prove now the second part. Let
$C_{1}$ ¨ $C_{2}$ be two chord diagrams of degrees
$k$ and $l$, respectively.
Consider the operators $O(C_{1})$ and $O(C_{2})$,
obtained be contracting the
$Tr(x_{k}\dots x_{1})$ and $Tr(x_{k+l},\dots,x_{k+1})$
at pairs of variables among
$x_{1},\dots, x_{k}$ and $x_{k+1},\dots,x_{k+l}$
according to diagrams $C_{1}$ and
$C_{2}$, respectively, see Fig.\ref{form}.
Taking the traces of these operators, we get
just  $U(C_{1})$ and $U(C_{2})$.

\begin{figure}
\begin{center}
\unitlength 1mm
\linethickness{0.4pt}
\begin{picture}(153.67,60.00)
\multiput(20.00,60.00)(1.07,-0.12){2}{\line(1,0){1.07}}
\multiput(22.15,59.77)(0.34,-0.12){6}{\line(1,0){0.34}}
\multiput(24.20,59.08)(0.19,-0.11){10}{\line(1,0){0.19}}
\multiput(26.05,57.96)(0.12,-0.11){13}{\line(1,0){0.12}}
\multiput(27.62,56.47)(0.11,-0.16){11}{\line(0,-1){0.16}}
\multiput(28.84,54.68)(0.11,-0.29){7}{\line(0,-1){0.29}}
\multiput(29.64,52.68)(0.12,-0.71){3}{\line(0,-1){0.71}}
\put(29.99,50.54){\line(0,-1){2.16}}
\multiput(29.87,48.38)(-0.12,-0.42){5}{\line(0,-1){0.42}}
\multiput(29.29,46.30)(-0.11,-0.21){9}{\line(0,-1){0.21}}
\multiput(28.28,44.39)(-0.12,-0.14){12}{\line(0,-1){0.14}}
\multiput(26.88,42.74)(-0.16,-0.12){11}{\line(-1,0){0.16}}
\multiput(25.16,41.43)(-0.25,-0.11){8}{\line(-1,0){0.25}}
\multiput(23.19,40.52)(-0.53,-0.12){4}{\line(-1,0){0.53}}
\put(21.08,40.06){\line(-1,0){2.16}}
\multiput(18.92,40.06)(-0.53,0.12){4}{\line(-1,0){0.53}}
\multiput(16.81,40.52)(-0.25,0.11){8}{\line(-1,0){0.25}}
\multiput(14.84,41.43)(-0.16,0.12){11}{\line(-1,0){0.16}}
\multiput(13.12,42.74)(-0.12,0.14){12}{\line(0,1){0.14}}
\multiput(11.72,44.39)(-0.11,0.21){9}{\line(0,1){0.21}}
\multiput(10.71,46.30)(-0.12,0.42){5}{\line(0,1){0.42}}
\put(10.13,48.38){\line(0,1){2.16}}
\multiput(10.01,50.54)(0.12,0.71){3}{\line(0,1){0.71}}
\multiput(10.36,52.68)(0.11,0.29){7}{\line(0,1){0.29}}
\multiput(11.16,54.68)(0.11,0.16){11}{\line(0,1){0.16}}
\multiput(12.38,56.47)(0.12,0.11){13}{\line(1,0){0.12}}
\multiput(13.95,57.96)(0.19,0.11){10}{\line(1,0){0.19}}
\multiput(15.80,59.08)(0.52,0.12){8}{\line(1,0){0.52}}
\multiput(60.00,60.00)(1.07,-0.12){2}{\line(1,0){1.07}}
\multiput(62.15,59.77)(0.34,-0.12){6}{\line(1,0){0.34}}
\multiput(64.20,59.08)(0.19,-0.11){10}{\line(1,0){0.19}}
\multiput(66.05,57.96)(0.12,-0.11){13}{\line(1,0){0.12}}
\multiput(67.62,56.47)(0.11,-0.16){11}{\line(0,-1){0.16}}
\multiput(68.84,54.68)(0.11,-0.29){7}{\line(0,-1){0.29}}
\multiput(69.64,52.68)(0.12,-0.71){3}{\line(0,-1){0.71}}
\put(69.99,50.54){\line(0,-1){2.16}}
\multiput(69.87,48.38)(-0.12,-0.42){5}{\line(0,-1){0.42}}
\multiput(69.29,46.30)(-0.11,-0.21){9}{\line(0,-1){0.21}}
\multiput(68.28,44.39)(-0.12,-0.14){12}{\line(0,-1){0.14}}
\multiput(66.88,42.74)(-0.16,-0.12){11}{\line(-1,0){0.16}}
\multiput(65.16,41.43)(-0.25,-0.11){8}{\line(-1,0){0.25}}
\multiput(63.19,40.52)(-0.53,-0.12){4}{\line(-1,0){0.53}}
\put(61.08,40.06){\line(-1,0){2.16}}
\multiput(58.92,40.06)(-0.53,0.12){4}{\line(-1,0){0.53}}
\multiput(56.81,40.52)(-0.25,0.11){8}{\line(-1,0){0.25}}
\multiput(54.84,41.43)(-0.16,0.12){11}{\line(-1,0){0.16}}
\multiput(53.12,42.74)(-0.12,0.14){12}{\line(0,1){0.14}}
\multiput(51.72,44.39)(-0.11,0.21){9}{\line(0,1){0.21}}
\multiput(50.71,46.30)(-0.12,0.42){5}{\line(0,1){0.42}}
\put(50.13,48.38){\line(0,1){2.16}}
\multiput(50.01,50.54)(0.12,0.71){3}{\line(0,1){0.71}}
\multiput(50.36,52.68)(0.11,0.29){7}{\line(0,1){0.29}}
\multiput(51.16,54.68)(0.11,0.16){11}{\line(0,1){0.16}}
\multiput(52.38,56.47)(0.12,0.11){13}{\line(1,0){0.12}}
\multiput(53.95,57.96)(0.19,0.11){10}{\line(1,0){0.19}}
\multiput(55.80,59.08)(0.52,0.12){8}{\line(1,0){0.52}}
\put(19.67,40.00){\line(0,1){20.00}}
\put(10.00,49.67){\line(1,0){20.00}}
\put(50.00,49.67){\line(1,0){19.67}}
\put(63.67,49.67){\line(1,0){6.33}}
\put(56.00,41.00){\line(0,1){18.33}}
\put(63.33,59.33){\line(0,-1){18.67}}
\put(29.00,45.67){\circle*{1.33}}
\put(68.67,45.33){\circle*{1.33}}

\put(-10.00,10.00){\line(1,0){30.00}}
\put(0.67,10.00){\oval(15.33,10.00)[t]}
\put(9.17,10.00){\oval(18.33,10.00)[t]}
\put(42.67,10.00){\oval(18.67,10.00)[t]}
\put(51.33,10.00){\oval(18.67,11.33)[t]}
\put(52.33,10.00){\oval(10.00,10.00)[t]}
\put(70.00,10.00){\line(1,0){30.00}}
\put(80.67,10.00){\oval(15.33,10.00)[t]}
\put(89.17,10.00){\oval(18.33,10.00)[t]}
\put(27.00,10.00){\line(1,0){39.00}}
\put(115.33,10.00){\oval(18.67,10.00)[t]}
\put(124.00,10.00){\oval(18.67,11.33)[t]}
\put(125.00,10.00){\oval(10.00,10.00)[t]}
\put(99.67,10.00){\line(1,0){39.00}}
\put(5.00,4.33){\makebox(0,0)[cc]{$O(C_1)$}}
\put(45.67,4.33){\makebox(0,0)[cc]{$O(C_2)$}}
\put(112.67,5.00){\makebox(0,0)[cc]{$O(C_{1}C_{2})$}}
\put(4.67,35.00){\makebox(0,0)[cc]{$C_{2}$}}
\put(45.00,35.00){\makebox(0,0)[cc]{$C_{2}$}}
\end{picture}

\end{center}
\label{form}
\caption{The Product Formula}
\end{figure}
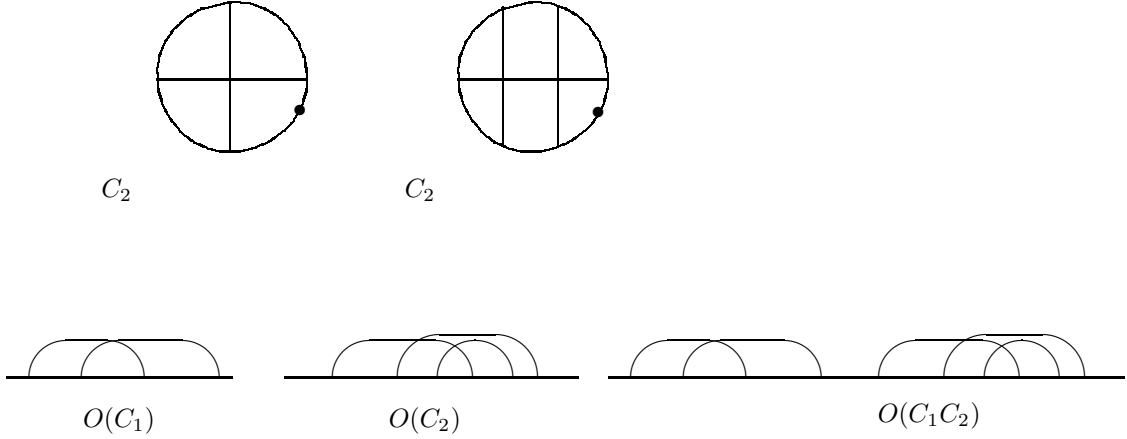

Naturally, these operators are adjoint invariant and
hence, scalar.
Thus we obtain

$$U(C_{1})\cdot U(C_{2})= Tr(O(C_{1})) \cdot Tr(O(C_{2}))=$$

$$ (n^{2}-1) Tr (O(C_{1}\cdot C_{2})) =
(n^{2}-1) U(C_{1}\cdot C_{1}), $$

That completes the proof of the theorem.
\end{proof}

\begin{exs}
Finally we give the list of chord diagrams of
orders $\le 4$ with values of $U$ on them.
Taking into account theorem 7, it is sufficient
to give polynomials on diagrams which are
not products, see fig. \ref{notpr}.

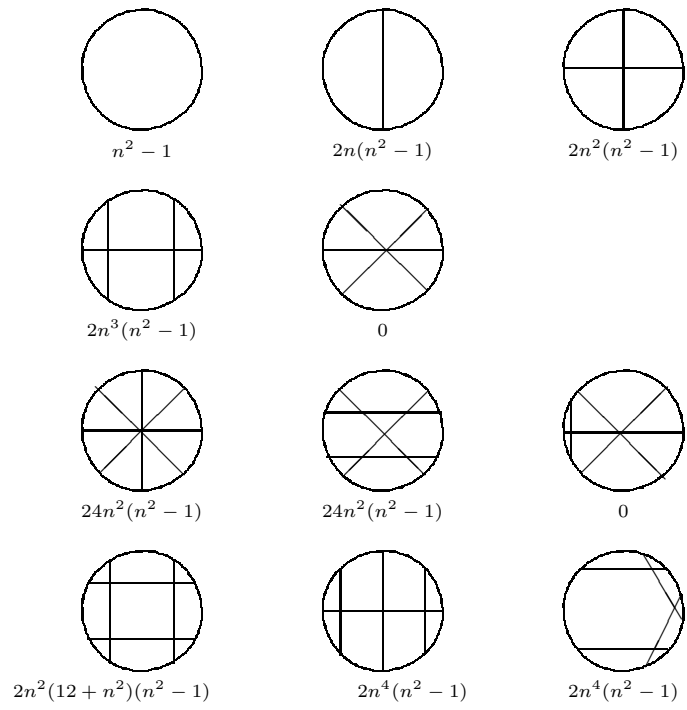
\begin{figure}
\unitlength 0.8mm
\linethickness{0.4pt}
\begin{center}
\begin{picture}(105.01,140.5)
\multiput(15,140.01)(1.08,-0.12){2}{\line(1,0){1.08}}
\multiput(17.15,139.77)(0.34,-0.12){6}{\line(1,0){0.34}}
\multiput(19.2,139.08)(0.19,-0.11){10}{\line(1,0){0.19}}
\multiput(21.05,137.97)(0.12,-0.11){13}{\line(1,0){0.12}}
\multiput(22.63,136.48)(0.11,-0.16){11}{\line(0,-1){0.16}}
\multiput(23.84,134.69)(0.11,-0.29){7}{\line(0,-1){0.29}}
\multiput(24.64,132.68)(0.12,-0.71){3}{\line(0,-1){0.71}}
\put(24.99,130.54){\line(0,-1){2.16}}
\multiput(24.87,128.38)(-0.12,-0.42){5}{\line(0,-1){0.42}}
\multiput(24.3,126.3)(-0.11,-0.21){9}{\line(0,-1){0.21}}
\multiput(23.28,124.39)(-0.12,-0.14){12}{\line(0,-1){0.14}}
\multiput(21.88,122.74)(-0.16,-0.12){11}{\line(-1,0){0.16}}
\multiput(20.16,121.43)(-0.25,-0.11){8}{\line(-1,0){0.25}}
\multiput(18.2,120.52)(-0.53,-0.12){4}{\line(-1,0){0.53}}
\put(16.08,120.05){\line(-1,0){2.16}}
\multiput(13.92,120.05)(-0.53,0.12){4}{\line(-1,0){0.53}}
\multiput(11.81,120.52)(-0.25,0.11){8}{\line(-1,0){0.25}}
\multiput(9.84,121.42)(-0.16,0.12){11}{\line(-1,0){0.16}}
\multiput(8.12,122.73)(-0.12,0.14){12}{\line(0,1){0.14}}
\multiput(6.72,124.38)(-0.11,0.21){9}{\line(0,1){0.21}}
\multiput(5.71,126.29)(-0.12,0.42){5}{\line(0,1){0.42}}
\put(5.13,128.38){\line(0,1){2.16}}
\multiput(5.01,130.54)(0.12,0.71){3}{\line(0,1){0.71}}
\multiput(5.36,132.67)(0.11,0.29){7}{\line(0,1){0.29}}
\multiput(6.16,134.68)(0.11,0.16){11}{\line(0,1){0.16}}
\multiput(7.37,136.47)(0.12,0.11){13}{\line(1,0){0.12}}
\multiput(8.94,137.96)(0.19,0.11){10}{\line(1,0){0.19}}
\multiput(10.79,139.08)(0.34,0.12){6}{\line(1,0){0.34}}
\multiput(12.84,139.77)(1.08,0.12){2}{\line(1,0){1.08}}
\multiput(55,140.01)(1.08,-0.12){2}{\line(1,0){1.08}}
\multiput(57.15,139.77)(0.34,-0.12){6}{\line(1,0){0.34}}
\multiput(59.2,139.08)(0.19,-0.11){10}{\line(1,0){0.19}}
\multiput(61.05,137.97)(0.12,-0.11){13}{\line(1,0){0.12}}
\multiput(62.63,136.48)(0.11,-0.16){11}{\line(0,-1){0.16}}
\multiput(63.84,134.69)(0.11,-0.29){7}{\line(0,-1){0.29}}
\multiput(64.64,132.68)(0.12,-0.71){3}{\line(0,-1){0.71}}
\put(64.99,130.54){\line(0,-1){2.16}}
\multiput(64.87,128.38)(-0.12,-0.42){5}{\line(0,-1){0.42}}
\multiput(64.3,126.3)(-0.11,-0.21){9}{\line(0,-1){0.21}}
\multiput(63.28,124.39)(-0.12,-0.14){12}{\line(0,-1){0.14}}
\multiput(61.88,122.74)(-0.16,-0.12){11}{\line(-1,0){0.16}}
\multiput(60.16,121.43)(-0.25,-0.11){8}{\line(-1,0){0.25}}
\multiput(58.2,120.52)(-0.53,-0.12){4}{\line(-1,0){0.53}}
\put(56.08,120.05){\line(-1,0){2.16}}
\multiput(53.92,120.05)(-0.53,0.12){4}{\line(-1,0){0.53}}
\multiput(51.81,120.52)(-0.25,0.11){8}{\line(-1,0){0.25}}
\multiput(49.84,121.42)(-0.16,0.12){11}{\line(-1,0){0.16}}
\multiput(48.12,122.73)(-0.12,0.14){12}{\line(0,1){0.14}}
\multiput(46.72,124.38)(-0.11,0.21){9}{\line(0,1){0.21}}
\multiput(45.71,126.29)(-0.12,0.42){5}{\line(0,1){0.42}}
\put(45.13,128.38){\line(0,1){2.16}}
\multiput(45.01,130.54)(0.12,0.71){3}{\line(0,1){0.71}}
\multiput(45.36,132.67)(0.11,0.29){7}{\line(0,1){0.29}}
\multiput(46.16,134.68)(0.11,0.16){11}{\line(0,1){0.16}}
\multiput(47.37,136.47)(0.12,0.11){13}{\line(1,0){0.12}}
\multiput(48.94,137.96)(0.19,0.11){10}{\line(1,0){0.19}}
\multiput(50.79,139.08)(0.34,0.12){6}{\line(1,0){0.34}}
\multiput(52.84,139.77)(1.08,0.12){2}{\line(1,0){1.08}}
\multiput(95,140.01)(1.08,-0.12){2}{\line(1,0){1.08}}
\multiput(97.15,139.77)(0.34,-0.12){6}{\line(1,0){0.34}}
\multiput(99.2,139.08)(0.19,-0.11){10}{\line(1,0){0.19}}
\multiput(101.05,137.97)(0.12,-0.11){13}{\line(1,0){0.12}}
\multiput(102.63,136.48)(0.11,-0.16){11}{\line(0,-1){0.16}}
\multiput(103.84,134.69)(0.11,-0.29){7}{\line(0,-1){0.29}}
\multiput(104.64,132.68)(0.12,-0.71){3}{\line(0,-1){0.71}}
\put(104.99,130.54){\line(0,-1){2.16}}
\multiput(104.87,128.38)(-0.12,-0.42){5}{\line(0,-1){0.42}}
\multiput(104.3,126.3)(-0.11,-0.21){9}{\line(0,-1){0.21}}
\multiput(103.28,124.39)(-0.12,-0.14){12}{\line(0,-1){0.14}}
\multiput(101.88,122.74)(-0.16,-0.12){11}{\line(-1,0){0.16}}
\multiput(100.16,121.43)(-0.25,-0.11){8}{\line(-1,0){0.25}}
\multiput(98.2,120.52)(-0.53,-0.12){4}{\line(-1,0){0.53}}
\put(96.08,120.05){\line(-1,0){2.16}}
\multiput(93.92,120.05)(-0.53,0.12){4}{\line(-1,0){0.53}}
\multiput(91.81,120.52)(-0.25,0.11){8}{\line(-1,0){0.25}}
\multiput(89.84,121.42)(-0.16,0.12){11}{\line(-1,0){0.16}}
\multiput(88.12,122.73)(-0.12,0.14){12}{\line(0,1){0.14}}
\multiput(86.72,124.38)(-0.11,0.21){9}{\line(0,1){0.21}}
\multiput(85.71,126.29)(-0.12,0.42){5}{\line(0,1){0.42}}
\put(85.13,128.38){\line(0,1){2.16}}
\multiput(85.01,130.54)(0.12,0.71){3}{\line(0,1){0.71}}
\multiput(85.36,132.67)(0.11,0.29){7}{\line(0,1){0.29}}
\multiput(86.16,134.68)(0.11,0.16){11}{\line(0,1){0.16}}
\multiput(87.37,136.47)(0.12,0.11){13}{\line(1,0){0.12}}
\multiput(88.94,137.96)(0.19,0.11){10}{\line(1,0){0.19}}
\multiput(90.79,139.08)(0.34,0.12){6}{\line(1,0){0.34}}
\multiput(92.84,139.77)(1.08,0.12){2}{\line(1,0){1.08}}
\multiput(55,110.01)(1.08,-0.12){2}{\line(1,0){1.08}}
\multiput(57.15,109.77)(0.34,-0.12){6}{\line(1,0){0.34}}
\multiput(59.2,109.08)(0.19,-0.11){10}{\line(1,0){0.19}}
\multiput(61.05,107.97)(0.12,-0.11){13}{\line(1,0){0.12}}
\multiput(62.63,106.48)(0.11,-0.16){11}{\line(0,-1){0.16}}
\multiput(63.84,104.69)(0.11,-0.29){7}{\line(0,-1){0.29}}
\multiput(64.64,102.68)(0.12,-0.71){3}{\line(0,-1){0.71}}
\put(64.99,100.54){\line(0,-1){2.16}}
\multiput(64.87,98.38)(-0.12,-0.42){5}{\line(0,-1){0.42}}
\multiput(64.3,96.3)(-0.11,-0.21){9}{\line(0,-1){0.21}}
\multiput(63.28,94.39)(-0.12,-0.14){12}{\line(0,-1){0.14}}
\multiput(61.88,92.74)(-0.16,-0.12){11}{\line(-1,0){0.16}}
\multiput(60.16,91.43)(-0.25,-0.11){8}{\line(-1,0){0.25}}
\multiput(58.2,90.52)(-0.53,-0.12){4}{\line(-1,0){0.53}}
\put(56.08,90.05){\line(-1,0){2.16}}
\multiput(53.92,90.05)(-0.53,0.12){4}{\line(-1,0){0.53}}
\multiput(51.81,90.52)(-0.25,0.11){8}{\line(-1,0){0.25}}
\multiput(49.84,91.42)(-0.16,0.12){11}{\line(-1,0){0.16}}
\multiput(48.12,92.73)(-0.12,0.14){12}{\line(0,1){0.14}}
\multiput(46.72,94.38)(-0.11,0.21){9}{\line(0,1){0.21}}
\multiput(45.71,96.29)(-0.12,0.42){5}{\line(0,1){0.42}}
\put(45.13,98.38){\line(0,1){2.16}}
\multiput(45.01,100.54)(0.12,0.71){3}{\line(0,1){0.71}}
\multiput(45.36,102.67)(0.11,0.29){7}{\line(0,1){0.29}}
\multiput(46.16,104.68)(0.11,0.16){11}{\line(0,1){0.16}}
\multiput(47.37,106.47)(0.12,0.11){13}{\line(1,0){0.12}}
\multiput(48.94,107.96)(0.19,0.11){10}{\line(1,0){0.19}}
\multiput(50.79,109.08)(0.34,0.12){6}{\line(1,0){0.34}}
\multiput(52.84,109.77)(1.08,0.12){2}{\line(1,0){1.08}}
\put(55,120.33){\line(0,1){19.67}}
\put(95,140){\line(0,-1){19.67}}
\put(85.33,130.33){\line(1,0){19.67}}
\multiput(15,110.01)(1.08,-0.12){2}{\line(1,0){1.08}}
\multiput(17.15,109.77)(0.34,-0.12){6}{\line(1,0){0.34}}
\multiput(19.2,109.08)(0.19,-0.11){10}{\line(1,0){0.19}}
\multiput(21.05,107.97)(0.12,-0.11){13}{\line(1,0){0.12}}
\multiput(22.63,106.48)(0.11,-0.16){11}{\line(0,-1){0.16}}
\multiput(23.84,104.69)(0.11,-0.29){7}{\line(0,-1){0.29}}
\multiput(24.64,102.68)(0.12,-0.71){3}{\line(0,-1){0.71}}
\put(24.99,100.54){\line(0,-1){2.16}}
\multiput(24.87,98.38)(-0.12,-0.42){5}{\line(0,-1){0.42}}
\multiput(24.3,96.3)(-0.11,-0.21){9}{\line(0,-1){0.21}}
\multiput(23.28,94.39)(-0.12,-0.14){12}{\line(0,-1){0.14}}
\multiput(21.88,92.74)(-0.16,-0.12){11}{\line(-1,0){0.16}}
\multiput(20.16,91.43)(-0.25,-0.11){8}{\line(-1,0){0.25}}
\multiput(18.2,90.52)(-0.53,-0.12){4}{\line(-1,0){0.53}}
\put(16.08,90.05){\line(-1,0){2.16}}
\multiput(13.92,90.05)(-0.53,0.12){4}{\line(-1,0){0.53}}
\multiput(11.81,90.52)(-0.25,0.11){8}{\line(-1,0){0.25}}
\multiput(9.84,91.42)(-0.16,0.12){11}{\line(-1,0){0.16}}
\multiput(8.12,92.73)(-0.12,0.14){12}{\line(0,1){0.14}}
\multiput(6.72,94.38)(-0.11,0.21){9}{\line(0,1){0.21}}
\multiput(5.71,96.29)(-0.12,0.42){5}{\line(0,1){0.42}}
\put(5.13,98.38){\line(0,1){2.16}}
\multiput(5.01,100.54)(0.12,0.71){3}{\line(0,1){0.71}}
\multiput(5.36,102.67)(0.11,0.29){7}{\line(0,1){0.29}}
\multiput(6.16,104.68)(0.11,0.16){11}{\line(0,1){0.16}}
\multiput(7.37,106.47)(0.12,0.11){13}{\line(1,0){0.12}}
\multiput(8.94,107.96)(0.19,0.11){10}{\line(1,0){0.19}}
\multiput(10.79,109.08)(0.34,0.12){6}{\line(1,0){0.34}}
\multiput(12.84,109.77)(1.08,0.12){2}{\line(1,0){1.08}}
\put(5,100){\line(1,0){20}}
\put(9.33,91.67){\line(0,1){16.67}}
\put(20.33,108.33){\line(0,-1){16.67}}
\put(48.33,92.67){\line(1,1){14}}
\put(62.33,93.33){\line(-1,1){14.33}}
\put(45,100){\line(1,0){20}}
\multiput(15,80.01)(1.08,-0.12){2}{\line(1,0){1.08}}
\multiput(17.15,79.77)(0.34,-0.12){6}{\line(1,0){0.34}}
\multiput(19.2,79.08)(0.19,-0.11){10}{\line(1,0){0.19}}
\multiput(21.05,77.97)(0.12,-0.11){13}{\line(1,0){0.12}}
\multiput(22.63,76.48)(0.11,-0.16){11}{\line(0,-1){0.16}}
\multiput(23.84,74.69)(0.11,-0.29){7}{\line(0,-1){0.29}}
\multiput(24.64,72.68)(0.12,-0.71){3}{\line(0,-1){0.71}}
\put(24.99,70.54){\line(0,-1){2.16}}
\multiput(24.87,68.38)(-0.12,-0.42){5}{\line(0,-1){0.42}}
\multiput(24.3,66.3)(-0.11,-0.21){9}{\line(0,-1){0.21}}
\multiput(23.28,64.39)(-0.12,-0.14){12}{\line(0,-1){0.14}}
\multiput(21.88,62.74)(-0.16,-0.12){11}{\line(-1,0){0.16}}
\multiput(20.16,61.43)(-0.25,-0.11){8}{\line(-1,0){0.25}}
\multiput(18.2,60.52)(-0.53,-0.12){4}{\line(-1,0){0.53}}
\put(16.08,60.05){\line(-1,0){2.16}}
\multiput(13.92,60.05)(-0.53,0.12){4}{\line(-1,0){0.53}}
\multiput(11.81,60.52)(-0.25,0.11){8}{\line(-1,0){0.25}}
\multiput(9.84,61.42)(-0.16,0.12){11}{\line(-1,0){0.16}}
\multiput(8.12,62.73)(-0.12,0.14){12}{\line(0,1){0.14}}
\multiput(6.72,64.38)(-0.11,0.21){9}{\line(0,1){0.21}}
\multiput(5.71,66.29)(-0.12,0.42){5}{\line(0,1){0.42}}
\put(5.13,68.38){\line(0,1){2.16}}
\multiput(5.01,70.54)(0.12,0.71){3}{\line(0,1){0.71}}
\multiput(5.36,72.67)(0.11,0.29){7}{\line(0,1){0.29}}
\multiput(6.16,74.68)(0.11,0.16){11}{\line(0,1){0.16}}
\multiput(7.37,76.47)(0.12,0.11){13}{\line(1,0){0.12}}
\multiput(8.94,77.96)(0.19,0.11){10}{\line(1,0){0.19}}
\multiput(10.79,79.08)(0.34,0.12){6}{\line(1,0){0.34}}
\multiput(12.84,79.77)(1.08,0.12){2}{\line(1,0){1.08}}
\multiput(55,80.01)(1.08,-0.12){2}{\line(1,0){1.08}}
\multiput(57.15,79.77)(0.34,-0.12){6}{\line(1,0){0.34}}
\multiput(59.2,79.08)(0.19,-0.11){10}{\line(1,0){0.19}}
\multiput(61.05,77.97)(0.12,-0.11){13}{\line(1,0){0.12}}
\multiput(62.63,76.48)(0.11,-0.16){11}{\line(0,-1){0.16}}
\multiput(63.84,74.69)(0.11,-0.29){7}{\line(0,-1){0.29}}
\multiput(64.64,72.68)(0.12,-0.71){3}{\line(0,-1){0.71}}
\put(64.99,70.54){\line(0,-1){2.16}}
\multiput(64.87,68.38)(-0.12,-0.42){5}{\line(0,-1){0.42}}
\multiput(64.3,66.3)(-0.11,-0.21){9}{\line(0,-1){0.21}}
\multiput(63.28,64.39)(-0.12,-0.14){12}{\line(0,-1){0.14}}
\multiput(61.88,62.74)(-0.16,-0.12){11}{\line(-1,0){0.16}}
\multiput(60.16,61.43)(-0.25,-0.11){8}{\line(-1,0){0.25}}
\multiput(58.2,60.52)(-0.53,-0.12){4}{\line(-1,0){0.53}}
\put(56.08,60.05){\line(-1,0){2.16}}
\multiput(53.92,60.05)(-0.53,0.12){4}{\line(-1,0){0.53}}
\multiput(51.81,60.52)(-0.25,0.11){8}{\line(-1,0){0.25}}
\multiput(49.84,61.42)(-0.16,0.12){11}{\line(-1,0){0.16}}
\multiput(48.12,62.73)(-0.12,0.14){12}{\line(0,1){0.14}}
\multiput(46.72,64.38)(-0.11,0.21){9}{\line(0,1){0.21}}
\multiput(45.71,66.29)(-0.12,0.42){5}{\line(0,1){0.42}}
\put(45.13,68.38){\line(0,1){2.16}}
\multiput(45.01,70.54)(0.12,0.71){3}{\line(0,1){0.71}}
\multiput(45.36,72.67)(0.11,0.29){7}{\line(0,1){0.29}}
\multiput(46.16,74.68)(0.11,0.16){11}{\line(0,1){0.16}}
\multiput(47.37,76.47)(0.12,0.11){13}{\line(1,0){0.12}}
\multiput(48.94,77.96)(0.19,0.11){10}{\line(1,0){0.19}}
\multiput(50.79,79.08)(0.34,0.12){6}{\line(1,0){0.34}}
\multiput(52.84,79.77)(1.08,0.12){2}{\line(1,0){1.08}}
\multiput(95,80.01)(1.08,-0.12){2}{\line(1,0){1.08}}
\multiput(97.15,79.77)(0.34,-0.12){6}{\line(1,0){0.34}}
\multiput(99.2,79.08)(0.19,-0.11){10}{\line(1,0){0.19}}
\multiput(101.05,77.97)(0.12,-0.11){13}{\line(1,0){0.12}}
\multiput(102.63,76.48)(0.11,-0.16){11}{\line(0,-1){0.16}}
\multiput(103.84,74.69)(0.11,-0.29){7}{\line(0,-1){0.29}}
\multiput(104.64,72.68)(0.12,-0.71){3}{\line(0,-1){0.71}}
\put(104.99,70.54){\line(0,-1){2.16}}
\multiput(104.87,68.38)(-0.12,-0.42){5}{\line(0,-1){0.42}}
\multiput(104.3,66.3)(-0.11,-0.21){9}{\line(0,-1){0.21}}
\multiput(103.28,64.39)(-0.12,-0.14){12}{\line(0,-1){0.14}}
\multiput(101.88,62.74)(-0.16,-0.12){11}{\line(-1,0){0.16}}
\multiput(100.16,61.43)(-0.25,-0.11){8}{\line(-1,0){0.25}}
\multiput(98.2,60.52)(-0.53,-0.12){4}{\line(-1,0){0.53}}
\put(96.08,60.05){\line(-1,0){2.16}}
\multiput(93.92,60.05)(-0.53,0.12){4}{\line(-1,0){0.53}}
\multiput(91.81,60.52)(-0.25,0.11){8}{\line(-1,0){0.25}}
\multiput(89.84,61.42)(-0.16,0.12){11}{\line(-1,0){0.16}}
\multiput(88.12,62.73)(-0.12,0.14){12}{\line(0,1){0.14}}
\multiput(86.72,64.38)(-0.11,0.21){9}{\line(0,1){0.21}}
\multiput(85.71,66.29)(-0.12,0.42){5}{\line(0,1){0.42}}
\put(85.13,68.38){\line(0,1){2.16}}
\multiput(85.01,70.54)(0.12,0.71){3}{\line(0,1){0.71}}
\multiput(85.36,72.67)(0.11,0.29){7}{\line(0,1){0.29}}
\multiput(86.16,74.68)(0.11,0.16){11}{\line(0,1){0.16}}
\multiput(87.37,76.47)(0.12,0.11){13}{\line(1,0){0.12}}
\multiput(88.94,77.96)(0.19,0.11){10}{\line(1,0){0.19}}
\multiput(90.79,79.08)(0.34,0.12){6}{\line(1,0){0.34}}
\multiput(92.84,79.77)(1.08,0.12){2}{\line(1,0){1.08}}
\multiput(15,50.01)(1.08,-0.12){2}{\line(1,0){1.08}}
\multiput(17.15,49.77)(0.34,-0.12){6}{\line(1,0){0.34}}
\multiput(19.2,49.08)(0.19,-0.11){10}{\line(1,0){0.19}}
\multiput(21.05,47.97)(0.12,-0.11){13}{\line(1,0){0.12}}
\multiput(22.63,46.48)(0.11,-0.16){11}{\line(0,-1){0.16}}
\multiput(23.84,44.69)(0.11,-0.29){7}{\line(0,-1){0.29}}
\multiput(24.64,42.68)(0.12,-0.71){3}{\line(0,-1){0.71}}
\put(24.99,40.54){\line(0,-1){2.16}}
\multiput(24.87,38.38)(-0.12,-0.42){5}{\line(0,-1){0.42}}
\multiput(24.3,36.3)(-0.11,-0.21){9}{\line(0,-1){0.21}}
\multiput(23.28,34.39)(-0.12,-0.14){12}{\line(0,-1){0.14}}
\multiput(21.88,32.74)(-0.16,-0.12){11}{\line(-1,0){0.16}}
\multiput(20.16,31.43)(-0.25,-0.11){8}{\line(-1,0){0.25}}
\multiput(18.2,30.52)(-0.53,-0.12){4}{\line(-1,0){0.53}}
\put(16.08,30.05){\line(-1,0){2.16}}
\multiput(13.92,30.05)(-0.53,0.12){4}{\line(-1,0){0.53}}
\multiput(11.81,30.52)(-0.25,0.11){8}{\line(-1,0){0.25}}
\multiput(9.84,31.42)(-0.16,0.12){11}{\line(-1,0){0.16}}
\multiput(8.12,32.73)(-0.12,0.14){12}{\line(0,1){0.14}}
\multiput(6.72,34.38)(-0.11,0.21){9}{\line(0,1){0.21}}
\multiput(5.71,36.29)(-0.12,0.42){5}{\line(0,1){0.42}}
\put(5.13,38.38){\line(0,1){2.16}}
\multiput(5.01,40.54)(0.12,0.71){3}{\line(0,1){0.71}}
\multiput(5.36,42.67)(0.11,0.29){7}{\line(0,1){0.29}}
\multiput(6.16,44.68)(0.11,0.16){11}{\line(0,1){0.16}}
\multiput(7.37,46.47)(0.12,0.11){13}{\line(1,0){0.12}}
\multiput(8.94,47.96)(0.19,0.11){10}{\line(1,0){0.19}}
\multiput(10.79,49.08)(0.34,0.12){6}{\line(1,0){0.34}}
\multiput(12.84,49.77)(1.08,0.12){2}{\line(1,0){1.08}}
\multiput(55,50.01)(1.08,-0.12){2}{\line(1,0){1.08}}
\multiput(57.15,49.77)(0.34,-0.12){6}{\line(1,0){0.34}}
\multiput(59.2,49.08)(0.19,-0.11){10}{\line(1,0){0.19}}
\multiput(61.05,47.97)(0.12,-0.11){13}{\line(1,0){0.12}}
\multiput(62.63,46.48)(0.11,-0.16){11}{\line(0,-1){0.16}}
\multiput(63.84,44.69)(0.11,-0.29){7}{\line(0,-1){0.29}}
\multiput(64.64,42.68)(0.12,-0.71){3}{\line(0,-1){0.71}}
\put(64.99,40.54){\line(0,-1){2.16}}
\multiput(64.87,38.38)(-0.12,-0.42){5}{\line(0,-1){0.42}}
\multiput(64.3,36.3)(-0.11,-0.21){9}{\line(0,-1){0.21}}
\multiput(63.28,34.39)(-0.12,-0.14){12}{\line(0,-1){0.14}}
\multiput(61.88,32.74)(-0.16,-0.12){11}{\line(-1,0){0.16}}
\multiput(60.16,31.43)(-0.25,-0.11){8}{\line(-1,0){0.25}}
\multiput(58.2,30.52)(-0.53,-0.12){4}{\line(-1,0){0.53}}
\put(56.08,30.05){\line(-1,0){2.16}}
\multiput(53.92,30.05)(-0.53,0.12){4}{\line(-1,0){0.53}}
\multiput(51.81,30.52)(-0.25,0.11){8}{\line(-1,0){0.25}}
\multiput(49.84,31.42)(-0.16,0.12){11}{\line(-1,0){0.16}}
\multiput(48.12,32.73)(-0.12,0.14){12}{\line(0,1){0.14}}
\multiput(46.72,34.38)(-0.11,0.21){9}{\line(0,1){0.21}}
\multiput(45.71,36.29)(-0.12,0.42){5}{\line(0,1){0.42}}
\put(45.13,38.38){\line(0,1){2.16}}
\multiput(45.01,40.54)(0.12,0.71){3}{\line(0,1){0.71}}
\multiput(45.36,42.67)(0.11,0.29){7}{\line(0,1){0.29}}
\multiput(46.16,44.68)(0.11,0.16){11}{\line(0,1){0.16}}
\multiput(47.37,46.47)(0.12,0.11){13}{\line(1,0){0.12}}
\multiput(48.94,47.96)(0.19,0.11){10}{\line(1,0){0.19}}
\multiput(50.79,49.08)(0.34,0.12){6}{\line(1,0){0.34}}
\multiput(52.84,49.77)(1.08,0.12){2}{\line(1,0){1.08}}
\multiput(95,50.01)(1.08,-0.12){2}{\line(1,0){1.08}}
\multiput(97.15,49.77)(0.34,-0.12){6}{\line(1,0){0.34}}
\multiput(99.2,49.08)(0.19,-0.11){10}{\line(1,0){0.19}}
\multiput(101.05,47.97)(0.12,-0.11){13}{\line(1,0){0.12}}
\multiput(102.63,46.48)(0.11,-0.16){11}{\line(0,-1){0.16}}
\multiput(103.84,44.69)(0.11,-0.29){7}{\line(0,-1){0.29}}
\multiput(104.64,42.68)(0.12,-0.71){3}{\line(0,-1){0.71}}
\put(104.99,40.54){\line(0,-1){2.16}}
\multiput(104.87,38.38)(-0.12,-0.42){5}{\line(0,-1){0.42}}
\multiput(104.3,36.3)(-0.11,-0.21){9}{\line(0,-1){0.21}}
\multiput(103.28,34.39)(-0.12,-0.14){12}{\line(0,-1){0.14}}
\multiput(101.88,32.74)(-0.16,-0.12){11}{\line(-1,0){0.16}}
\multiput(100.16,31.43)(-0.25,-0.11){8}{\line(-1,0){0.25}}
\multiput(98.2,30.52)(-0.53,-0.12){4}{\line(-1,0){0.53}}
\put(96.08,30.05){\line(-1,0){2.16}}
\multiput(93.92,30.05)(-0.53,0.12){4}{\line(-1,0){0.53}}
\multiput(91.81,30.52)(-0.25,0.11){8}{\line(-1,0){0.25}}
\multiput(89.84,31.42)(-0.16,0.12){11}{\line(-1,0){0.16}}
\multiput(88.12,32.73)(-0.12,0.14){12}{\line(0,1){0.14}}
\multiput(86.72,34.38)(-0.11,0.21){9}{\line(0,1){0.21}}
\multiput(85.71,36.29)(-0.12,0.42){5}{\line(0,1){0.42}}
\put(85.13,38.38){\line(0,1){2.16}}
\multiput(85.01,40.54)(0.12,0.71){3}{\line(0,1){0.71}}
\multiput(85.36,42.67)(0.11,0.29){7}{\line(0,1){0.29}}
\multiput(86.16,44.68)(0.11,0.16){11}{\line(0,1){0.16}}
\multiput(87.37,46.47)(0.12,0.11){13}{\line(1,0){0.12}}
\multiput(88.94,47.96)(0.19,0.11){10}{\line(1,0){0.19}}
\multiput(90.79,49.08)(0.34,0.12){6}{\line(1,0){0.34}}
\multiput(92.84,49.77)(1.08,0.12){2}{\line(1,0){1.08}}
\put(5,70){\line(1,0){20}}
\put(15,60){\line(0,1){20}}
\put(22,77){\line(-1,-1){14}}
\put(22,62.67){\line(-1,1){14.67}}
\put(47.67,77){\line(1,-1){14.33}}
\put(62.67,76.67){\line(-1,-1){14.33}}
\put(45.33,73){\line(1,0){19.33}}
\put(64.33,65.67){\line(-1,0){18.67}}
\put(88,63){\line(1,1){14.33}}
\put(87.33,76.67){\line(1,-1){14.67}}
\put(104.67,69.67){\line(-1,0){19.67}}
\put(86.33,65){\line(0,1){10.33}}
\put(9.67,31.33){\line(0,1){17}}
\put(20.33,48.33){\line(0,-1){17}}
\put(23.67,35.33){\line(-1,0){17.67}}
\put(6,44.67){\line(1,0){17.67}}
\put(45,40){\line(1,0){20}}
\put(48,32.67){\line(0,1){14}}
\put(55,50){\line(0,-1){20}}
\put(62,33){\line(0,1){14}}
\put(87.67,47){\line(1,0){14.67}}
\put(98.33,49.33){\line(3,-5){6.4}}
\put(104.67,42.67){\line(-1,-2){6}}
\put(102.33,33.67){\line(-1,0){15}}
\put(15,116.67){\makebox(0,0)[cc]{\scriptsize{$n^{2}-1$}}}
\put(55,116.67){\makebox(0,0)[cc]{\scriptsize{$2n(n^{2}-1)$}}}
\put(95,116.67){\makebox(0,0)[cc]{\scriptsize{$2n^{2}(n^{2}-1)$}}}
\put(15,86.67){\makebox(0,0)[cc]{\scriptsize{$2n^{3}(n^{2}-1)$}}}
\put(15,56.67){\makebox(0,0)[cc]{\scriptsize{$24n^{2}(n^{2}-1)$}}}
\put(10,26.67){\makebox(0,0)[cc]{\scriptsize{$2n^{2}(12+n^{2})(n^{2}-1)$}}}
\put(55,86.67){\makebox(0,0)[cc]{\scriptsize{$0$}}}
\put(55,56.67){\makebox(0,0)[cc]{\scriptsize{$24n^{2}(n^{2}-1)$}}}
\put(60,26.67){\makebox(0,0)[cc]{\scriptsize{$2n^{4}(n^{2}-1)$}}}
\put(95,56.67){\makebox(0,0)[cc]{\scriptsize{$0$}}}
\put(95,26.67){\makebox(0,0)[cc]{\scriptsize{$2n^{4}(n^{2}-1)$}}}

\end{picture}
\end{center}

\caption{Values $U$ on chord diagrams}
\label{notpr}

\end{figure}

\end{exs}

\newpage

\section{The case of $sp(n)$}

Like in the preceding cases, to determine $(3)$
we have to compute the trace formulae $(4)$ and $(5)$.

\begin{lm}
For $sp(n)$ the sum (4) looks like
$${1\over 2} (Tr P Tr Q- Tr (PjQ^{*}j))\eqno(12),$$
and (5) looks like $${1\over 2}(Tr (PQ)+ Tr (PjQ^{*}j)),\eqno(13)$$

where $j$ is the block--diagonal $2n\times 2n$--matrix:

$$j=\left(\begin{array}{cc} 0 & E \\ -E & 0 \end{array}\right). $$
\end{lm}

\begin{proof}
Consider the basis given by the matrices $E_{ii}-E_{n+i,n+i}\;\left(  1\leq
i\leq n\right),$  $E_{ij}-E_{n+j,n+i}\;\left(  i\neq j\right)$,

$E_{i,n+i},E_{n+i,i}\;\left(  1\leq i\leq n\right)$  $E_{i,n+j}+E_{j,n+i},$
$E_{n+i,j}+E_{n+j,i}\;\left(  1\leq i<j\leq n\right)  $.

Respect to
it, the dual basis is given by $\frac{1}{2}\left(  E_{ii}-E_{n+i,n+i}\right)
,\;\frac{1}{2}\left(  E_{ji}-E_{n+i,n+j}\right),$ $E_{n+i,i},\;E_{i,n+i}%
,\;\frac{1}{2}\left(  E_{n+i,j}+E_{n+j,i}\right)$, \newline $\frac{1}{2}\left(
E_{i,n+j}+E_{j,n+i}\right)  $. We now determine the trace formulae by
evaluating them on distinct pairs $\left(  P,Q\right)  :$

\begin{enumerate}
\item $$\left(  P,Q\right)  =\left(  id_{2n},id_{2n}\right)  $$

$$\sum_{\alpha}Tr\left(  \alpha\alpha^{\prime}\right)  =\dim\;sp\left(
n\right)  =2n^{2}+n$$

$$\sum_{\alpha}Tr\left(  \alpha\right)  Tr\left(
\alpha^{\prime}\right)  =0$$

\item $$\left(  P,Q\right)  =\left(  j,j\right)  $$

$$\sum_{\alpha
}Tr\left(  j\alpha j\alpha^{\prime}\right)  =n$$

$$\sum_{a}Tr\left(
j\alpha\right)  Tr\left(  j\alpha^{\prime}\right)  =2\sum_{i=1}^{n}Tr\left(
-E_{ii}\right)  Tr\left(  E_{n+i,n+i}\right)  =-2n$$

\item
$$\left(  P,Q\right)  =\left(  E_{11},E_{11}\right)  $$

$$\sum
_{a}Tr\left(  E_{11}jE_{11}\alpha^{\prime}\right)  =\sum_{a}Tr\left(
E_{11}\alpha\right)  Tr\left(  E_{11}\alpha^{\prime}\right)  =\frac{1}{2}$$
\end{enumerate}

Taking into account the obvious trace formulae for
the identity matrix and the
matrix $j$, the assertion follows at once.
\end{proof}

Substituing in equation (3) we obtain:
\begin{eqnarray*}
f(x_{1},\dots, x_{k})=
{1\over 2}\sum_{{l=0} \atop {i_{1}>\dots >i_{l},i_{l+1}<\dots<i_{n}}} ^{n}
(-1)^{l}\times
(Tr (x_{i_1}\dots x_{i_{k}}) Tr (x_{i_{k+1}}\dots x_{i_n})+\\
Tr(x_{i_1}\dots x_{i_{k}} x_{i_n}^{*}\dots x_{i_{k+1}}^{*}))
\end{eqnarray*}

\begin{th}
Each $m$--linear form on $x_{i_1},\dots x_{i_m}$, obtained from
the adjoint representation for $sp(n)$ is a
linear combination of traces of products for $x_{i_j}$ and $x_{i_j}^{*}$
with polynomial coefficients with respect to $n$.
\end{th}

\makebox[1in][l]{R.Campoamor-Stursberg}\newline
\makebox[1in][l]{Departamento de Geometr\'{\i}a y Topolog\'{\i}a}\newline
\makebox[1in][l]{Facultad CC. Matem\'aticas}\newline
\makebox[1in][l]{Plaza de Ciencias, 3}\newline
\makebox[1in][l]{Universidad Complutense}\newline
\makebox[1in][l]{E-28040 Madrid (Spain)}\newline
\makebox[1in][l]{rutwig{@}mat.ucm.es}\newline

\makebox[1in][l]{V.O.Manturov}\newline
\makebox[1in][l]{Moscow State University}\newline
\makebox[1in][l]{Faculty of Mechanics and Mathematics}\newline
\makebox[1in][l]{Vorobyovy gory}\newline
\makebox[1in][l]{Ru-119899 Moscow)}\newline
\makebox[1in][l]{vassily{@}manturov.mccme.ru}

\end{document}